\input pictex.tex
\input amstex.tex
\documentstyle{amsppt}
\magnification=1200
 \pagewidth{16.3truecm}
 \pageheight{24truecm}
 \nologo

\def\sing{\operatorname{sing}}

\def\dim{\operatorname{dim}}

\def\reg{\operatorname{reg}}
\def\red{\operatorname{red}}

\refstyle{A}
\widestnumber\key{ACLM2}
\topmatter
\title
Arc spaces, motivic integration and stringy invariants
\endtitle
\author
Willem Veys  
\endauthor
\address K.U.Leuven, Departement Wiskunde, Celestijnenlaan 200B,
         B--3001 Leuven, Belgium  \endaddress
\email wim.veys\@wis.kuleuven.ac.be  \newline
 http://www.wis.kuleuven.ac.be/algebra/veys.htm
\endemail
\endtopmatter
\document

\bigskip
The concept of {\it motivic integration} was invented by
Kontsevich to show that birationally equivalent Calabi-Yau
manifolds have the same Hodge numbers. He constructed a certain
measure on the {\it arc space} of an algebraic variety, the
motivic measure, with the subtle and crucial property that it
takes values not in $\Bbb R$, but in the Grothendieck ring of
algebraic varieties. A whole theory on this subject was then
developed by Denef and Loeser in various papers, with several
applications.

Batyrev introduced with motivic integration techniques  new
singularity invariants, the {\it stringy invariants}, for
algebraic varieties with mild singularities, more precisely log
terminal singularities. He used them for instance to formulate a
topological Mirror Symmetry test for pairs of singular Calabi-Yau
varieties. We generalized these invariants to almost arbitrary
singular varieties, assuming Mori's Minimal Model Program.

The aim of these notes is to provide a gentle introduction to
these concepts.
 There exist already good surveys by Denef-Loeser
[DL8] and Looijenga [Loo], and a nice elementary introduction by
Craw [Cr]. Here we merely want to explain the basic concepts and
first results, including the $p$-adic number theoretic pre-history
of the theory, and to provide concrete examples.

The text is a slightly adapted version of the \lq extended
abstract\rq\ of the author's talks at the 12th MSJ-IRI
"Singularity Theory and Its Applications" (2003) in Sapporo. At
the end we included a list of various recent results.

\bigskip
\bigskip
\noindent {\bf Table of contents }
\bigskip
\noindent 1. Pre-history  \dotfill 1 \vskip1pt

\noindent 2. Arc spaces \dotfill 3 \vskip1pt

\noindent 3. Motivic integration \dotfill 8 \vskip1pt

\noindent 4. First applications \dotfill 12 \vskip1pt

\noindent 5. Motivic volume \dotfill 13 \vskip1pt

\noindent 6. Motivic zeta functions \dotfill 15 \vskip1pt

\noindent 7. Batyrev's stringy invariants \dotfill 18 \vskip1pt

\noindent 8. Stringy invariants for general singularities \dotfill
22 \vskip1pt

\noindent 9. Miscellaneous recent results \dotfill 27 \vskip1pt

\noindent References \dotfill 29 \vskip1pt

\pagebreak

\noindent {\bf 1. Pre-history}
\bigskip
\noindent {\bf 1.1.} Let $f \in \Bbb Z[x_1,\cdots,x_m]$ and $r \in
\Bbb Z_{> 0}$. A very general problem in number theory is to
compute the number of solutions of the congruence
$f(x_1,\cdots,x_m) = 0$ mod $r$ (in $(\Bbb Z/r \Bbb Z)^m$). Thanks
to the Chinese remainder theorem it is enough to consider the case
where $r$ is a power of a prime.

So we fix a prime number $p$ and we investigate congruences modulo
varying powers of $p$. We denote by $F_n$ the number of solutions
of $f(x_1,\cdots,x_m) = 0$ mod $p^{n+1}$.
\bigskip
\noindent {\bf 1.2.} {\sl Examples.}
\medskip
\noindent \item{1.} $f_1 = y - x^2$. It should be clear that $F_n
= p^{n+1}$.
 \medskip \item{2.} $f_2 = x \cdot y$. {\smc Exercise} : $F_n
= (n+2)p^{n+1} - (n+1) p^n$. \medskip \item{3.} $f_3 = y^2 - x^3$.
We list $F_n$ for small $n$ : $F_0 = p$,
\medskip
\settabs 19 \columns \+ $F_1 = p(2p - 1)$ & & & & $F_5 = p^5 (p^2
+ p - 1)$ & & & & & $F_7 = p^7 (2p^2 - 1)$  & & & & & $F_{11} =
p^{11} (p^3 + p^2 - 1)$ \cr \+  $F_2 = p^2(2p-1)$ & & & & $F_6 =
p^6(p^2 + p - 1)$ & & & & & $F_8 = p^8(2 p^2 - 1)$ & & & & &
$F_{12} = p^{12} (p^3 + p^2 - 1)$. \cr \+  $F_3 = p^3 (2p-1)$ & &
& & & & & & & $F_9 = p^9 (2 p^2 - 1)$ \cr \+  $F_4 = p^4 (2p - 1)$
& & & & & & & & & $F_{10} = p^{10} (2 p^2 - 1)$ \cr

\bigskip
Note that the plane curve $\{ f_1  = 0 \}$ is nonsingular, $\{ f_2
= 0 \}$ has the easiest curve singularity, an ordinary node, and
$\{ f_3 = 0 \}$ has a slightly more complicated singularity, an
ordinary cusp.  It is in fact this cusp which is responsible for
the at first sight not so nice behavior of the $F_n$ for $f_3$.

More generally, the problem of the behavior of the $F_n$ turns out
to be non-obvious precisely when $\{ f = 0 \}$ has singularities.

\bigskip
\noindent {\bf 1.3.} We now know that, for any $f \in \Bbb Z[x_1,
\cdots, x_m ]$, the $F_n$ do satisfy the following `regular'
behavior.

\proclaim {Conjecture {\rm[Borewicz, Shafarevich]} = Theorem {\rm[Igusa]}}
The generating formal series $J_p (T) := J_p(f,T) = \sum_{n \geq
0} F_n T^n$ is a rational function in $T$. (In particular the
$F_n$ are determined by a {\rm finite} number of them.)
\endproclaim

\smallskip
\noindent Igusa showed this in 1975 [Ig1] using

\smallskip
(1) a `translation' of $J_p(T)$ into a $p$-adic integral (more
precisely into $\int_{\Bbb Z^m_p} |f|^s_p |dx|$, which is now
called {\it Igusa's local zeta function}, and which is the
ancestor of the motivic zeta function of section 6),

(2) an embedded
resolution of singularities for $\{ f = 0 \}$,

(3) the
change of variables formula for integrals.

\noindent (We will see later an analogue of this strategy in the
theory of motivic integration.)
\bigskip
\noindent {\bf 1.4.} {\sl Examples} (continuing 1.2).

\smallskip
\item{1.} $J_p(f_1;T) = \frac{p}{1 - pT}$ (easy). \item{2.} {\smc
Exercise} : $J_p(f_2;T) = \frac{2p-1-p^2T}{(1 - pT)^2}$. \item{3.}
Claim : $J_p(f_3;T) = p \frac{1 + (p-1) T + (p^6 - p^5)T^5 - p^7
T^6}{(1 - p^7 T^6)(1 - pT)}$.

\bigskip
\noindent {\bf 1.5.} We already want to mention another connection
with singularity theory; the famous (still open) {\it monodromy
conjecture} of Igusa relates the poles of $J_p(T)$ with
eigenvalues of local monodromy of $f$ considered as a map $f :
\Bbb C^n \rightarrow \Bbb C$, see (6.8).

\bigskip
\noindent {\bf 1.6.} Before introducing arc spaces and motivic
integration in the next sections, we present a hopefully motivating
analogy between this number theoretic setting and the geometric
arc setting.

\bigskip
\settabs 2 \columns \+ $f \in \Bbb Z[x_1, \cdots , x_m]$ & $f \in
\Bbb C[x_1, \cdots , x_m]$ \cr \medskip \noindent \hrule \medskip
\+ solution of $f = 0$ over the ring & solution of $f = 0$ over
the ring \cr \+ $\Bbb Z/p^{n+1} \Bbb Z \cong \Bbb Z_p/p^{n+1} \Bbb
Z_p$, i.e. & $\Bbb C[t]/(t^{n+1}) \cong \Bbb C[[t]]/(t^{n+1})$,
i.e. \cr \+ an $m$-tuple with coordinates of the form & an
$m$-tuple with coordinates of the form \cr \+ $a_0 + a_1 p + ... +
a_n p^n \   (a_i \in \{ 0,1,...,p-1 \})$ & $a_0 + a_1 t + ... +
a_n t^n \ \ (a_i \in \Bbb C)$ \cr \+ & (``n-jet" of $\{ f = 0 \}$)
\cr
\medskip
\noindent \hrule
\medskip
\+ solution of $f = 0$ over $\Bbb Z_p = \underset \leftarrow \to
\lim\, \Bbb Z/p^{n+1} \Bbb Z$, & solution of $f = 0$ over $\Bbb
C[[t]] = \underset \leftarrow \to \lim \,\Bbb C[t]/(t^{n+1}),$ \cr
\+ i.e. with coordinates of the form & i.e. with coordinates of
the form \cr \+ $\sum^\infty_{n=0} a_i p^i$ & $\sum^\infty_{n=0}
a_i t^i$ \cr \+ & (``arc" of $\{ f = 0 \}$) \cr
\medskip \noindent \hrule
\medskip \+ integrate over $\Bbb Z^m_p$ & integrate over ${\Cal
L}(\Bbb C^m) := \{ \text{ arcs of } \Bbb C^m \}$ \cr

\bigskip
\noindent {\sl Warning.} Here and further on we sometimes use
other (better ?) normalizations than in the original papers.

\bigskip
\bigskip
\noindent {\bf 2. Arc spaces}
\bigskip
\noindent Let $X$ be an algebraic variety over $\Bbb C$. (The
theory can be generalized to any field of characteristic zero.)
\bigskip
\noindent {\bf 2.1.} The {\it space of arcs modulo} $t^{n+1}$ or {\it
space of $n$-jets on} $X$ is an algebraic variety ${\Cal L}_n
(X)$ over $\Bbb C$ such that
$$\{\text{points of } {\Cal L}_n(X) \text { with coordinates in }
\Bbb C \} = \{\text {points of } X \text{ with coordinates in }
\frac{\Bbb C[t]}{(t^{n+1})} \}.$$
For all $n$ there are obvious
`truncation maps' $\pi^{n+1}_n : {\Cal L}_{n+1} (X) \rightarrow
{\Cal L}_n (X)$, obtained by reducing $(n+1)$-jets modulo
$t^{n+1}$, and more generally $\pi^m_n : {\Cal L}_m(X) \rightarrow
{\Cal L}_n(X)$ for $m \geq n$. This description is somewhat
informal, but is essentially what is needed. We now first provide
examples and give the `exact' definition later.
\bigskip
\noindent {\bf  2.2.} {\sl Example.} Let $X = \Bbb C^d$. Then
$$\split {\Cal L}_n (X) & = \{ (a^{(1)}_0 + a^{(1)}_1 t + \cdots + a^{(1)}_n
t^n, \cdots , a^{(d)}_0 + a^{(d)}_1 t + \cdots + a^{(d)}_n t^n),
\text { with all }  a^{(j)}_i \in \Bbb C \} \\ & \cong \Bbb
C^{(n+1)d}. \endsplit$$
\medskip
\noindent {\bf 2.3.} {\sl Example.} Let $X = \{ y^2 - x^3 = 0 \}$

\smallskip
\noindent
(0) ${\Cal L}_0 (X) = \{ (a_0,b_0) \in \Bbb C^2 | b^2_0 =
a^3_0 \} = X$.
\smallskip
\noindent
(1) $$\align {\Cal L}_1 (X) & = \{ (a_0 +
a_1 t, b_0 + b_1 t) \in (\Bbb C[t]/(t^2))^2 \mid (b_0 + b_1 t)^2 = (a_0 +
a_1 t)^3 \text{mod } t^2 \} \\ & = \{ (a_0 + a_1 t, b_0 + b_1 t)
\in (\Bbb C[t]/(t^2))^2 \mid b^2_0 = a^3_0 \text{ and } 2b_0 b_1 = 3a^2_0
a_1 \}.
\endalign$$
So we can consider ${\Cal L}_1(X)$ as the (2-dimensional)
algebraic variety in $\Bbb C^4$ with equations $b^2_0 = a^3_0$ and
$2b_0 b_1 = 3 a^2_0 a_1$ in the coordinates $a_0, a_1, b_0, b_1$.
The map $\pi^1_0 : {\Cal L}_1 (X) \rightarrow {\Cal L}_0 (X) = X$
is induced by the projection $\Bbb C^4 \rightarrow \Bbb C^2 :
(a_0, a_1, b_0, b_1) \mapsto (a_0, b_0)$.

 The fibre of $\pi_0^1$
above $(0,0)$ is $\{ (0,a_1,0,b_1) \} \cong \Bbb C^2$; this
corresponds to the fact that the tangent space to $X$ at $(0,0)$
is the whole $\Bbb C^2$. The fibre above $(a_0,b_0) \ne (0,0)$ is
the line in the $(a_1,b_1)$-plane with equation $2b_0 b_1 = 3
a^2_0 a_1$, which corresponds to the tangent line at $X$ in
$(a_0,b_0)$. In other words : ${\Cal L}_1(X)$ is the tangent
bundle $TX$, and $\pi_0^1$ is the natural projection $TX
\rightarrow X$.

\smallskip
\noindent
(2) ${\Cal L}_2(X) = \{ (a_0 + a_1 t + a_2 t^2, b_0 + b_1 t
+ b_2 t^2) \in (\Bbb C[t]/(t^3))^2  \mid (b_0 + b_1 t + b_2 t^2)^2 = (a_0
+ a_1 t + a_2 t^2)^3 \text {mod } t^3 \}$ is given in $\Bbb C^6$
by the equations
$$\left\{ \aligned & b^2_0  = a^3_0 \\
& 2b_0 b_1  = 3a^2_0 a_1 \\
& b^2_1 + 2b_0b_2 = 3a_0 a^2_1 + 3a^2_0 a_2.
\endaligned
\right.
$$

\noindent {\smc Exercise}. a) Verify the description of ${\Cal
L}_2 (X)$ and note that the map $\pi_1^2 : {\Cal L}_2
(X)\rightarrow {\Cal L}_1(X)$ is not surjective. More precisely,
the fibre of $\pi_0^2$ above $(0,0)$ is $\{ (0,a_1,a_2,0,0,b_2) \}
\cong \Bbb C^3$, but its image by $\pi_1^2$ is not the whole
$(a_1,b_1)$-plane; it is just the line $\{b_1=0\}$.

b) Compute ${\Cal L}_3(X)$ and note that also $\pi^3_2 : {\Cal
L}_3 (X) \rightarrow {\Cal L}_2(X)$ is not surjective.

c) However, above the nonsingular part of $X = {\Cal L}_0(X)$ all
considered maps $\pi^{n+1}_{n} : {\Cal L}_{n+1}(X) \rightarrow
{\Cal L}_n(X)$ are fibrations with fibre $\Bbb C$.

\bigskip
\noindent {\bf 2.4.} Some observations in the examples are easily
seen to be satisfied in general.

\noindent
(1) ${\Cal L}_0(X) = X,
\quad {\Cal L}_1(X) = TX.$

\noindent (2) If $X$ is smooth of
dimension $d$, then all $\pi^{n+1}_n$ are locally trivial
fibrations (w.r.t. the Zariski topology) with fibre $\Bbb C^d$.

\bigskip
\noindent {\bf 2.5.} The {\it space of arcs on} $X$ is an
`algebraic variety of infinite dimension' ${\Cal L}(X)$ over $\Bbb
C$ such that
$$\{ \text{points of } {\Cal L}(X) \text{ with coordinates in } \Bbb
C \} = \{ \text{points of } X \text{with coordinates in } \Bbb
C[[t]] \}.$$ We provide the `exact' definition after continuing
the examples.
Now we have for all $n$ truncation maps $\pi_n : {\Cal L}(X)
\rightarrow {\Cal L}_n(X)$, obtained by reducing arcs modulo
$t^{n+1}$.
\bigskip
\noindent {\bf 2.6.} {\sl Example.} Let $X = \Bbb C^d$. Then
$${\Cal L}(X) = \{ (\sum^\infty_{n=0} a^{(1)}_n t^n , \cdots ,
\sum^\infty_{n=0} a^{(d)}_n t^n), \text{ with all } a^{(j)}_n \in
\Bbb C \},$$ which can be considered as an infinite dimensional
affine space.

\bigskip
\noindent {\bf  2.7.} {\sl Example.} Let $X = \{ y^2 - x^3 = 0
\}$. Then ${\Cal L}(X)$ is given in the infinite dimensional
affine space with coordinates
$$\left\{ \aligned & a_0, a_1, a_2, \cdots , a_n, \cdots \\
& b_0, b_1, b_2, \cdots , b_n, \cdots \endaligned \right.
$$
by the infinite number of equations
$$\left\{ \aligned & b^2_0 = a^3_0 \\ & 2 b_0 b_1 = 3 a^2_0 a_1 \\
& b^2_1 + 2b_0 b_2 = 3a_0 a^2_1 + 3a^2_0 a_2 \\ & \cdots
\endaligned \right.$$
\bigskip
\noindent {\bf 2.8.} {\sl More precise definitions.}

\smallskip
(i) The `base extension operation' $Y \rightarrow Y \times_{\Bbb
C} \Bbb C[t]/(t^{n+1})$ is a covariant functor on the category of
complex algebraic varieties, and it has a right adjoint $X
\rightarrow {\Cal L}_n(X)$. This says that, for any $\Bbb
C$-algebra $R$, the set of $R$-valued points of ${\Cal L}_n(X)$ is
in natural bijection with the set of $R[t]/(t^{n+1})$-valued
points of $X$. In particular, as we said in (2.1), the $\Bbb
C$-valued points of ${\Cal L}_n(X)$ can be naturally identified
with the $\Bbb C[t]/(t^{n+1})$-valued points of $X$.

(ii) Then ${\Cal L}(X)$ is the inverse limit $\underset \leftarrow
\to \lim\, {\Cal L}_n(X)$. (Technically, it is important here that
the truncation morphisms $\pi^{n+1}_n : {\Cal L}_{n+1} (X)
\rightarrow {\Cal L}_n (X)$ are affine.)  The $K$-valued points of
${\Cal L}(X)$, for any field $K \supset \Bbb C$, are in natural
bijection with the $K[[t]]$-valued points of $X$. We mention the
following result, attributed to Kolchin : if $X$ is irreducible,
then ${\Cal L}(X)$ is irreducible.

See [DL3] for more information.
\bigskip
\noindent {\bf 2.9.} When $X$ is an affine variety, i.e. given by
a finite number of polynomial equations, one can describe
equations for the ${\Cal L}_n(X)$ and for ${\Cal L}(X)$ as in
Examples 2.3 and 2.7.
\bigskip
\noindent {\bf  2.10.} Some first natural and fundamental
questions are how the ${\Cal L}_n(X)$ and $\pi_n({\Cal L}(X))$
change with $n$. (For $\pi_n({\Cal L}(X))$ this was already
considered by Nash [Na].) Note that ${\Cal L}_n(X)$ describes by
definition the $n$-jets on $X$, and $\pi_n({\Cal L}(X))$ those
$n$-jets that can be lifted to arcs on $X$.

This can be compared with the number theoretical setting of the
previous section : there the question was how the solutions over
$\Bbb Z/p^{n+1} \Bbb Z$ changed with $n$, and we could consider
the same question for those solutions over $\Bbb Z/p^{n+1} \Bbb Z$
that can be lifted to solutions over $\Bbb Z_p$.
\bigskip
\noindent {\bf 2.11.} We now introduce the Grothendieck ring of
algebraic varieties, which is the `best' framework to answer these
questions, and which is moreover (essentially) the value ring for
motivic integration, to be explained in the next section.

Recall first two fundamental properties of the topological Euler
characteristic $\chi(\cdot) \in \Bbb Z$ on complex algebraic
varieties :

(1) $\chi(V) = \chi (Z) + \chi (V \setminus Z)$
if $Z$ is (Zariski-)closed in $V$,

(2) $\chi (V \times  W)
= \chi(V) \cdot \chi (W)$.

\medskip \noindent A finer invariant
satisfying these properties is the Hodge-Deligne polynomial
$H(\cdot) = H(\cdot;u,v) \in \Bbb Z[u,v]$, given for an algebraic
variety $V$ of dimension $d$ by $$H(V;u,v) := \sum^d_{p,q=0}
(\sum^{2d}_{i=0} (-1)^i h^{p,q} (H^i_c(V,\Bbb C)))u^pv^q,$$ where
$h^{p,q}(\cdot)$ denotes the dimension of the $(p,q)$-component of
the mixed Hodge structure.  Note that $H(V;1,1) = \chi(V)$.

The Grothendieck ring is the value ring of the `universal Euler
characteristic' on algebraic varieties.
\bigskip
\noindent {\bf Definition.} (i) The Grothendieck group of
(complex) algebraic varieties is the abelian group $K_0(Var_{\Bbb
C})$ generated by symbols $[V]$, where $V$ is an algebraic
variety, with the relations $[V] = [W]$ if $V$ and $W$ are
isomorphic, and $[V] = [Z] + [V \setminus Z]$ if $Z$ is
(Zariski-) closed in $V$.

(ii) there is a natural ring structure on $K_0(Var_{\Bbb C})$
given by $[V] \cdot [W] := [V \times W]$.

\medskip
\noindent --- So by construction the map $\{\text{Varieties over }
\Bbb C \} \rightarrow K_0(Var_{\Bbb C}) : V \mapsto [V]$ is indeed
universal with respect to the two properties above. Of course we
still loose some information by this operation. For example $X =
\{ y^2 - x^3 = 0 \} \subset \Bbb A^2$ satisfies $[X] = [\Bbb
A^1]$. Also, when $V\to B$ is a locally trivial fibration with
fibre $F$, then $[V]=[B] \cdot [F]$.  ---
\medskip
(iii) Let $C$ be a constructible subset of some variety $V$, i.e.
a disjoint union of (finitely many) locally closed subvarieties
$A_i$ of $V$, then $[C] \in K_0(Var_{\Bbb C})$ is well defined as
$[C] := \sum_i [A_i]$.

(iv) We denote  1:= [point], $\Bbb L := [\Bbb A^1]$ and ${\Cal
M}_{\Bbb C} := K_0 (Var_{\Bbb C})_{\Bbb L}$ the ring obtained from
$K_0(Var_{\Bbb C})$ by inverting ${\Bbb L}$.

\bigskip The rings $K_0(Var_{\Bbb C})$ and ${\Cal M}_{\Bbb C}$ are
quite mysterious. For instance, it was shown only recently that
$K_0(Var_{\Bbb C})$ is not a domain [Po], and it is still not
known whether ${\Cal M}_{\Bbb C}$ is a domain or not, or whether
the natural map $K_0 (Var_{\Bbb C}) \rightarrow {\Cal M}_{\Bbb C}$
is injective.

\bigskip
\noindent {\sl Remark.} There is an interesting alternative
description of $K_0(Var_{\Bbb C})$ as the abelian group, generated
by isomorphism classes $[V]$ of {\sl nonsingular projective}
varieties $V$, with the relations $[\emptyset]=0$ and $[\tilde V]
- [E]=[V]-[Z]$, where $\tilde V \to V$ is the blowing-up with
centre $Z$ and exceptional variety $E$ [Bi1].

\bigskip \noindent {\bf 2.12.}
We now answer the questions in (2.10). We will consider $[{\Cal
L}_n(X)]$ and $[\pi_n({\Cal L}(X))]$ in ${\Cal M}_{\Bbb C}$. For
the latter we use a theorem of Greenberg [Gr], stating that
$\pi_n({\Cal L}(X))$ is a constructible subset of ${\Cal L}_n(X)$.

\bigskip
\noindent \proclaim{Theorem {\rm[DL3][DL8]}} The generating formal series
$$J(T) := \sum_{n \geq 0} [{\Cal L}_n(X)]T^n \text { and } P(T) :=
\sum_{n \geq 0} [ \pi_n ({\Cal L}(X))]T^n$$ in ${\Cal M}_{\Bbb C}
[[T]]$ are rational, with moreover as denominators products of
polynomials of the form $1 - \Bbb L^a T^b$, where $a \in \Bbb Z$
and $b \in \Bbb Z_{> 0}$. \endproclaim

\bigskip
The proof uses motivic integration, which \lq explains\rq\ why
${\Cal M}_{\Bbb C}$ is needed instead of $K_0(Var_{\Bbb C})$; see
section 3.

 This result specializes
to the analogous statement, replacing $[\cdot]$ by $\chi(\cdot)$
or $H(\cdot)$. Note for this that $\chi : K_0(Var_{\Bbb C})
\rightarrow \Bbb Z$ and $H : K_0(Var_{\Bbb C}) \rightarrow \Bbb
Z[u,v]$ obviously extend to $\chi : {\Cal M}_{\Bbb C} \rightarrow
\Bbb Z$ and $H : {\Cal M}_{\Bbb C} \rightarrow \Bbb Z[u,v][\frac
1{uv}]$.  When $X = \{ f = 0 \}$ for some polynomial $f$, the
statement for $J(T)$ should be compared with Theorem 1.3 for
$J_p(T)$ ! In this case, we will outline a proof for $J(T)$ later.
We just mention that the proof for $P(T)$ uses techniques from
logic, more precisely quantifier elimination.
\bigskip
\noindent {\bf 2.13.} {\sl Example.} When $X$ is smooth of
dimension $d$, all ${\Cal L}_n(X) = \pi_n ({\Cal L}(X))$ are
locally trivial over $X$ with fibre $\Bbb C^{nd}$. Hence
$$J(T) = P(T) = \sum_{n \geq 0} [X] \Bbb L^{nd} T^n = \frac{[X]}{1
- \Bbb L^d T}.$$

\bigskip
\noindent {\bf 2.14.} {\sl Example.} Let $X = \{ y^2 - x^3 = 0 \}$.
The descriptions in Example 2.3 yield $[ {\Cal L}_0(X)] = [X] =
\Bbb L, [{\Cal L}_1(X)] = \Bbb L^2 + (\Bbb L - 1) \Bbb L = 2 \Bbb
L^2 - \Bbb L, [{\Cal L}_2(X)] = \Bbb L^3 + (\Bbb L - 1) \Bbb L^2 =
2 \Bbb L^3 - \Bbb L^2$.  We claim that
$$J(T) = \Bbb L \frac{1 + (\Bbb L - 1)T + (\Bbb L^6 - \Bbb L^5)T^5
- \Bbb L^7 T^6}{(1 - \Bbb L^7 T^6)(1 - \Bbb L T)},$$ see section
6. (Compare with 1.4(3)!)
 The formula in [DL5, Proposition 10.2.1] yields
$$P(T) = \frac{\Bbb L  + (1 - \Bbb L)T
- \Bbb L T^2}{(1 - \Bbb L T)(1 - T^2)}.$$
\bigskip
\noindent {\bf  2.15.} {\sl Example.} Let $X = \{ xy = 0 \}$.
{\smc Exercise} :

(i) $[{\Cal L}_n(X)] = (n+2)\Bbb L^{n+1} - (n+1) \Bbb L^n$. Then
$$J(T) = \frac{2 \Bbb L - 1 - \Bbb L^2 T}{(1 - \Bbb L T)^2}.$$
(Compare again with Examples 1.2 and 1.4.)

(ii) $[\pi_n({\Cal L}(X))] = 2 \Bbb L^{n+1} - 1$. Then
$$P(T) = \frac{2 \Bbb L - 1 - \Bbb L T}{(1 - \Bbb L T)(1 - T)}.$$

\bigskip
\noindent {\bf 2.16.} [Mu1] To conclude this section, we relate
some properties of the spaces of $n$-jets on $X$ to properties of
$X$. Let $d$ denote the dimension of $X$.

\smallskip
\noindent (i) The closure in $\Cal L_n(X)$ of
$(\pi^n_0)^{-1}(X_{\reg})$ is an irreducible component of $\Cal
L_n(X)$ of dimension $d(n+1)$.

\smallskip
\noindent (ii) Suppose that $X$ is locally a complete
intersection. Then

(1) $\Cal L_n(X)$ is pure dimensional if and only if $\dim \Cal
L_n(X) \leq d(n+1)$.

(2) $\Cal L_n(X)$ is irreducible if and only if $\dim
(\pi^n_0)^{-1}(X_{\sing}) < d(n+1)$.

(3) If $\Cal L_{n+1}(X)$ is pure dimensional or irreducible, then
so is $\Cal L_n(X)$.

(4) If $\Cal L_n(X)$ is irreducible for some $n>0$, then $X$ is
normal.

(5) $\Cal L_n(X)$ is irreducible for all $n>0$ if and only if $X$
has rational singularities.

\smallskip
\noindent (iii) When $d=1$ we have for any $n>0$ that $\Cal
L_n(X)$ is irreducible if and only if $X$ is nonsingular.

\bigskip
\bigskip
\noindent {\bf 3. Motivic integration}
\bigskip
\noindent This
notion is due to Kontsevich [Ko] on nonsingular varieties. It has
been further developed by Batyrev [Ba2][Ba3], and especially by Denef
and Loeser [DL3][DL4][DL6][DL8], with some improvements by Looijenga [Loo].
Probably the best way to view and understand it, is as being an
analogue of $p$-adic integration.
\bigskip
\noindent Let in this section $X$ be any algebraic variety of pure
dimension $d$.
\bigskip
\noindent {\bf 3.1.} A subset $A$ of ${\Cal L}(X)$ is called {\it
constructible} or {\it cylindric} or {\it a cylinder} if $A =
\pi^{-1}_m C$ for some $m$ and some constructible subset $C$ of
${\Cal L}_m(X)$. These can be considered as `reasonably nice'
subsets of the arc space ${\Cal L}(X)$, being precisely all arcs
obtained by lifting a nice subset of a jet space.

\bigskip
\noindent {\bf 3.2.} Suppose that $X$ is nonsingular. Then such a
constructible subset $A = \pi^{-1}_m C$ satisfies the property
$$[\pi_n(A)] = \Bbb L^{(n-m)d} [C] \quad \text { for all } n \geq
m ,$$
since $\pi^n_m : {\Cal L}_n (X) = \pi_n ({\Cal L}(X))
\rightarrow {\Cal L}_m (X) = \pi_m({\Cal L}(X))$ is a locally
trivial fibration with fibre $\Bbb C^{(n-m)d}$. We have in
particular that the
$$\frac{[\pi_n(A)]}{\Bbb L^{nd}}$$
are all equal in ${\Cal M}_{\Bbb C}$ for $n \geq m$.

For general $X$, a constructible set $A \subset {\Cal L}(X)$ for
which $A \cap {\Cal L}(X_{sing}) = \emptyset$ still satisfies
the property that the $\frac{[\pi_n(A)]}{\Bbb L^{nd}}$ stabilize
for $n$ big enough [DL3, Lemma 4.1]. More precisely we have the following.
\bigskip
\noindent {\bf Definition.} We call a set $A \subset {\Cal L}(X)$
{\it stable} if for some $m \in \Bbb N$ we have

(i) $\pi_m(A)$ is constructible and $A = \pi^{-1}_m (\pi_m(A))$,
and

(ii) for all $n \geq m$ the projection $\pi_{n+1}(A) \rightarrow
\pi_n(A)$ is a piecewise trivial fibration with fiber $\Bbb C^d$.

\noindent (So in particular $A$ is constructible.)

\bigskip
\noindent \proclaim{Lemma {\rm [DL3]}} If $A \subset {\Cal L}(X)$ is
constructible and $A \cap {\Cal L}(X_{\text{sing}}) = \emptyset$,
then $A$ is stable. \endproclaim

\medskip
Hence for such $A$ it makes sense to consider $\lim_{n \rightarrow
\infty} \frac{[\pi_n(A)]}{\Bbb L^{nd}} \in {\Cal M}_{\Bbb C}$ as
an invariant of $A$; it is called its {\it naive motivic measure}.
Note that for nonsingular $X$ the measure of ${\Cal L}(X)$ is just
$[X]$.
\bigskip
\noindent {\bf 3.3.} For arbitrary constructible $A \subset {\Cal
L}(X)$ the sequence $\frac{[\pi_n(A)]}{\Bbb L^{nd}}$ will not
stabilize.
\bigskip
\noindent {\sl Example.} Let $X = \{ xy = 0 \}$. From Example 2.15
we see that
$$\frac{[\pi_n({\Cal L}(X))]}{\Bbb L^{nd}} = \frac{2 \Bbb
L^{n+1}-1}{\Bbb L^n} = 2 \Bbb L - \frac{1}{\Bbb L^n}.$$
This
sequence `almost' stabilizes (the singular point of $X$ of course
causes the trouble), and it would be nice to be able to consider
$2 \Bbb L$ as the limit of this sequence.

\bigskip
This will indeed work in Kontsevich's {\it completed Grothendieck
ring} $\hat {\Cal M}_{\Bbb C}$. This is by definition the
completion of ${\Cal M}_{\Bbb C}$ with respect to the decreasing
filtration $F^m, m \in \Bbb Z$, of ${\Cal M}_{\Bbb C}$, where
$F^m$ is the subgroup of ${\Cal M}_{\Bbb C}$ generated by the
elements $\frac{[S]}{\Bbb L^i}$ with $S$ an algebraic variety and
$\dim S - i \leq -m$. Note that this is indeed a ring filtration :
$F^m \cdot F^n \subset F^{m+n}$.
So $\hat {\Cal M}_{\Bbb C} = \underset \overleftarrow m \to \lim
\frac{{\Cal M}_{\Bbb C}}{F^m}$.

\noindent {\sl Continuing the example.} Indeed in $\hat {\Cal
M}_{\Bbb C}$ we have
$$\lim_{n \rightarrow \infty} \frac{[\pi_n({\Cal L}(X))]}{\Bbb
L^{nd}} = 2 \Bbb L - \lim_{n \rightarrow \infty} \frac{1}{\Bbb
L^n} = 2 \Bbb L.$$

\bigskip
\noindent \proclaim{Theorem {\rm [DL3]}} Let $A$ be a constructible
subset of $\Cal L(X)$. Then the limit
$$\mu(A) := \lim_{n \rightarrow \infty} \frac {[\pi_n(A)]}{\Bbb
L^{nd}}$$ exists in $\hat {\Cal M}_{\Bbb C}$.
\endproclaim
\bigskip
\noindent We call $\mu(A)$ the {\it motivic measure} of $A$. This
yields a $\sigma$-additive measure $\mu$ on the Boolean algebra of
constructible subsets of ${\Cal L}(X)$.

\bigskip \noindent {\sl Note.}
It is not known whether the natural map ${\Cal M}_{\Bbb C}
\rightarrow \hat {\Cal M}_{\Bbb C}$ is injective; its kernel is
$\cap_{m \in \Bbb Z} F^m$. However, e.g. the topological Euler
characteristic $\chi(\cdot)$ and the Hodge-Deligne polynomial
$H(\cdot)$ factor through the image of $\Cal M_{\Bbb C}$ in $\hat
{\Cal M}_{\Bbb C}$.

\bigskip\noindent
{\sl Remark.} Let $S \subsetneq X$ be a closed subvariety; it is
not difficult to see that $\Cal L (S)$ is {\sl not} a
constructible subset of $\Cal L (X)$. It is possible to introduce
more generally {\sl measurable} subsets of ${\Cal L}(X)$, and to
associate analogously a motivic measure (in $\hat {\Cal M}_{\Bbb
C}$) to those subsets [Ba2][DL6]; we then have that such $\Cal L
(S)$ are measurable of measure zero.



\bigskip
\noindent{\bf 3.4.} We briefly compare with the $p$-adic case. Let
$M$ be a $d$-dimensional submanifold of $\Bbb Z^m_p$, defined
algebraically. Denote by $|M(\Bbb Z / p^{n+1} \Bbb Z)|$ the
number of $\Bbb Z / p^{n+1}\Bbb Z$ (= $\Bbb Z_p / p^{n+1} \Bbb Z_p$)-valued points of $M$. Then $\frac {|M(\Bbb
Z/p^{n+1} \Bbb Z)|}{p^{(n+1)d}} \in \Bbb Z[ \frac 1 p ]$ is
constant for $n$ big enough and is called the volume $\mu_p(M)$ of
$M$.

For a singular $d$-dimensional subvariety $Z$ of $\Bbb Z^m_p$ one defines its {\it volume} as $\mu_p(Z) := \lim_{\epsilon \rightarrow 0} \mu_p (Z
\setminus T_\epsilon (Z_{\text{sing}})) \in \Bbb R$, where
$T_\epsilon$ denotes a small tubular neighbourhood \lq of radius $\epsilon$\rq. Then by a Theorem of
Oesterl\'e [Oe] we have
$$\mu_p (Z) = \lim_{n \rightarrow \infty} \frac{|Z(\Bbb Z/p^{n+1}
\Bbb Z)|}{p^{(n+1)d}} .$$
Note the analogy
\bigskip
\settabs 3 \columns \+ & $p$-adic & motivic \cr \medskip \hrule
\medskip
\+ integrate over & $\Bbb Z^m_p$ & $(\Bbb C[[t]])^m$ \cr \medskip
\noindent \hrule \medskip \+ value rings & $\Bbb Z$ &
$K_0(Var_{\Bbb C})$ \cr \+ & $\Bbb Z[ \frac 1 p ]$ & $\Cal M_{\Bbb
C}$ \cr \+ &  $\Bbb R$ & $\hat{\Cal M}_{\Bbb C}$ \cr
\bigskip
\noindent The brilliant idea of Kontsevich was to use $\hat {\Cal
M}_{\Bbb C}$ instead of $\Bbb R$ as a value ring for integration.

\bigskip
\noindent {\bf 3.5.} We can now consider in a natural way {\it
motivic integration}. We do not treat the most general setting;
the following suffices in practice. Let $A \subset {\Cal L}(X)$ be
constructible and $\alpha : A \rightarrow \Bbb Z \cup \{ + \infty
\}$ a function with constructible fibres $\alpha^{-1} \{ n \}, n
\in \Bbb Z$. Then
$$\int_A \Bbb L^{-\alpha} d \mu := \sum_{n \in \Bbb Z} \mu (\alpha^{-1} \{ n \}) \Bbb L^{-n}$$
in $\hat {\Cal M}_{\Bbb C}$,
whenever the right hand side converges in $\hat {\Cal M}_{\Bbb
C}$. Then we say that $\Bbb L^{-\alpha}$ is {\it integrable} on
$A$. (This will always be the case if $\alpha$ is bounded from
below.)
\bigskip
\noindent {\bf 3.6.} An important example of an integrable
function is induced by an effective Cartier divisor $D$ on $X$,
i.e. $D$ is an (eventually non-reduced) subvariety of $X$ which is
locally given by one equation.  Define $ord_t D : {\Cal L}(X)
\rightarrow \Bbb N \cup \{  + \infty \} : \gamma \mapsto ord_t
f_D(\gamma)$, where $f_D$ is a local equation of $D$ in a
neighbourhood of the origin $\pi_0(\gamma)$ of $\gamma$. Note e.g.
that $(ord_t D)(\gamma) = + \infty$ if and only if $\gamma \in
{\Cal L} (D_{\red})$ and $(ord_t D )( \gamma) = 0$ if and only if
$\pi_0 (\gamma) \not\in D_{\red}$. One easily verifies that $\Bbb
L^{-ord_tD}$ is integrable on ${\Cal L}(X)$.

We note that $(ord_t D)^{-1}(+\infty)=\Cal L(D_{\red})$ is {\sl not} constructible; it is however measurable with measure zero.

\bigskip
\noindent {\sl Example.} Take $X = \Bbb A^1$ and $D$ the divisor
associated to the function $x^N$, i.e. the `origin with
multiplicity $N$'.
\smallskip
\noindent {\smc Exercise}. (i) $N|(ord_t D)(\gamma)$ for all
$\gamma \in {\Cal L} (\Bbb A^1)$ and
$$\mu(\{ \gamma \in {\Cal L}(\Bbb A^1) \mid (ord_t D)(\gamma) = i N \}
) = \frac {\Bbb L-1}{\Bbb L^i} \text { for all } i \in \Bbb N.$$

(ii) $\int_{{\Cal L}(\Bbb A^1)} \Bbb L^{-ord_tD} d \mu =
\frac{(\Bbb L-1)\Bbb L^{N+1}}{\Bbb L^{N+1}-1} = (\Bbb L-1) + \frac
{ \Bbb L-1}{\Bbb L^{1+N}-1}$.

\smallskip\noindent
This example is the easiest case of the following very useful
formula.

\medskip
\noindent \proclaim{Proposition {\rm [Ba3][Cr]}} Let $X$ be
nonsingular and $D = \sum_{i \in S} N_i D_i$ a normal crossings
divisor on $X$, i.e. all $D_i$ are nonsingular hypersurfaces
intersecting transversely (and occurring with multiplicity $N_i$).
Denote $D^\circ_I := (\cap_{i \in I} D_i) \setminus (\cup_{\ell
\not \in I} D_\ell)$ for $I \subset S$; the $D^\circ_I, I \subset
S$, form a natural locally closed stratification of $X$ (note that
$D^\circ_\emptyset = X \setminus (\cup_{\ell \in S} D_\ell)$).
Then
$$\int_{{\Cal L(X)}} \Bbb L^{-ord_tD}d\mu = \sum_{I \subset S} [D^\circ_I]
\prod_{i \in I} \frac{\Bbb L-1}{\Bbb L^{1+N_i}-1}.$$
\endproclaim
\bigskip
\noindent {\bf 3.7.} The construction in (3.6) can be generalized
as follows. Let $\Cal I$ be a sheaf of ideals on $X$. Then we define
$$ord_t {\Cal I} : {\Cal L}(X) \rightarrow \Bbb N \cup \{ +
\infty \} : \gamma \mapsto \min_g ord_t g(\gamma),$$ where the
minimum is taken over $g \in {\Cal I}$ in a neighbourhood of
$\pi_0(\gamma)$. Of course, when ${\Cal I}$ is the ideal sheaf of
an effective Cartier divisor $D$, then $ord_t {\Cal I} = ord_t D$.
\bigskip
\noindent {\bf 3.8.} The most crucial ingredient in the theory of
motivic integration is the {\it change of variables formula} or
{\it transformation rule} for motivic integrals under a birational
morphism.

\bigskip
\noindent \proclaim{Theorem {\rm [DL3]}} (i) Let $h : Y \rightarrow X$ be a proper
birational morphism between algebraic varieties $X$ and $Y$, where
$Y$ is nonsingular. Let $A \subset {\Cal L}(X)$ be constructible
and $\alpha : A \rightarrow \Bbb Z \cup \{ + \infty \}$ such that
$\Bbb L^{-\alpha}$ is integrable on $A$. Then
$$\int_A \Bbb L^{-\alpha} d \mu = \int_{h^{-1}A} \Bbb L^{-(\alpha
\circ h) - ord_t(Jac_h)} d \mu .$$
Here the ideal sheaf $Jac_h$ is
defined as follows. When also $X$ is nonsingular, it is locally
generated by the `ordinary' Jacobian determinant with respect to
local coordinates on $X$ and $Y$. For general $X$, the sheaf of
regular differential $d$-forms $h^\ast(\Omega^d_X)$ is still a
submodule of $\Omega^d_Y$; but now $h^\ast (\Omega^d_X)$ is not
necessarily locally generated by one element. Taking (locally) a
generator $\omega_Y$ of $\Omega^d_Y$, each $h^\ast (\omega)$ for
$\omega \in \Omega^d_X$ can be written as $h^\ast (\omega) =
g_\omega \omega_Y$, and $Jac_h$ is defined as the ideal sheaf
which is (locally) generated by these $g_\omega$.

(ii) When also $X$ is nonsingular and $\alpha = ord_t D$ for some
effective divisor $D$ on $X$, we can rewrite the formula as
follows :
$$\int_A \Bbb L^{-ord_tD} d \mu = \int_{h^{-1}A} \Bbb L^{-ord_t(h^\ast D + K_{Y|X})}d \mu .$$
Here $h^\ast D$ is the pullback of $D$, i.e. locally given by the
equation $f \circ h$, if $D$ is given by the equation $f$. And
$K_{Y|X}$ is the relative canonical divisor, which is precisely
the effective divisor with equation the Jacobian determinant.
Alternatively, $K_{Y|X} = K_Y - h^\ast K_X$ where $K_\centerdot$
denotes the (ordinary) canonical divisor, i.e. the divisor of
zeros and poles of a differential $d$-form.
\endproclaim

\bigskip
\noindent {\sl Note.} The birational morphism $h$ above must be
proper in order to induce a bijection from $\Cal L(Y)$ to $\Cal
L(X)$ outside subsets of measure zero. More precisely, denoting by
$Exc$ the exceptional locus of $h$, we have a bijection from $\Cal
L(Y) \setminus \Cal L(Exc)$ to $\Cal L(X)\setminus \Cal
L(h(Exc))$. This is an easy consequence of the valuative criterion
of properness [Har, Theorem II.4.7].

\bigskip \noindent
{\smc Exercise}. Check the change of variables formula in the
following special case : $h$ is the blowing-up of a nonsingular
$X$ in a nonsingular centre, $A = {\Cal L}(X)$ and $\alpha$ is the
zero function.

\bigskip
\bigskip
\noindent {\bf 4. First applications}
\bigskip
\noindent {\bf 4.1.} Here we mean by a {\it Calabi-Yau manifold}
$M$ of dimension $d$ a nonsingular complete (=compact) algebraic variety,
which admits a nowhere vanishing regular differential $d$-form
$\omega_M$. Alternative formulations of this last condition are
that the first Chern class of the tangent bundle of $M$ is zero,
or that the canonical divisor $K_M$ of $M$ is zero.
\bigskip
\noindent \proclaim{Theorem {\rm[Ko]}} Let $X$ and $Y$ be birationally
equivalent Calabi-Yau manifolds. Then $[X] = [Y]$ in $\hat {\Cal
M}_{\Bbb C}$. \endproclaim
\medskip
\demo{Proof} Since $X$ and $Y$ are birationally equivalent
there exist a nonsingular complete algebraic variety $Z$ and
birational morphisms $h_X : Z \rightarrow X$ and $h_Y : Z
\rightarrow Y$. By the definition of the motivic measure and the
change of variables formula we have in $\hat {\Cal M}_{\Bbb C}$ :
$$[X] = \mu ({\Cal L}(X)) = \int_{{\Cal L}(X)} 1 d \mu =
\int_{{\Cal L}(Z)} \Bbb L^{-ord_t K_{Z|X}}d \mu = \int_{{\Cal
L}(Z)} \Bbb L^{-ord_t K_Z} d \mu ,$$
and of course [Y] is given by
the same right hand side.  \qed
\enddemo

\bigskip This implies that
birationally equivalent Calabi-Yau manifolds have the same
Hodge-Deligne polynomial, meaning that they have the same Hodge
numbers. This result was Kontsevich's motivation to invent motivic
integration !
\bigskip
\noindent The same proof gives the following more general result.
Two nonsingular complete algebraic varieties are called {\it
K-equivalent} if there exists a nonsingular complete algebraic
variety $Z$ and birational morphisms $h_X : Z \rightarrow X$ and
$h_Y : Z \rightarrow Y$ such that $h^\ast_X K_X = h^\ast_Y K_Y$.
This is an important notion in birational geometry.

\bigskip
\noindent \proclaim {Theorem} Let $X$ and $Y$ be $K$-equivalent
varieties. Then $[X] = [Y]$ in $\hat {\Cal M}_{\Bbb C}$.
\endproclaim

\bigskip \noindent {\bf 4.2.} Let $h : Y \rightarrow X$ be a proper
birational morphism between nonsingular algebraic varieties. We
assume that the exceptional locus $Exc$ of $h$, i.e. the
subvariety of $Y$ where $h$ is not an isomorphism, is a normal
crossings divisor. Let $E_i, i \in S$, be the irreducible
components of $Exc$. The relative canonical divisor $K_{Y|X}$ is
supported on $Exc$; let $\nu_i - 1$ be the multiplicity of $E_i$
in this divisor, so $K_{Y|X} = \sum_{i \in S} (\nu_i - 1)E_i$.
Denoting $E^\circ_I := (\cap_{i \in I} E_i) \setminus (\cup_{\ell
\not \in I} E_\ell)$ for $I \subset S$, we have

$$[X] = \sum_{I \subset S} [E^\circ_I] \prod_{i \in I}
\frac{\Bbb L-1}{\Bbb L^{\nu_i} - 1} = \sum_{I \subset S}
[E^\circ_I] \prod_{i \in I} \frac{1}{[ \Bbb P^{\nu_i - 1}]}$$
in
$\hat {\Cal M}_{\Bbb C}$. Indeed, by the change of variables
formula we have again that
$$[X] = \mu ({\Cal L}(X)) = \int_{{\Cal L}(Y)} \Bbb L^{-ord_t
K_{Y|X}}d \mu,$$
and then Proposition 3.6 yields the stated
formula. Specializing to the topological Euler characteristic
yields the remarkable formula
$$\chi(X) = \sum_{I \subset S} \chi(E^\circ_I) \prod_{i \in I}
\frac{1}{\nu_i} ,$$ which was first surprisingly obtained in [DL1],
using $p$-adic integration and the Grothen- dieck-Lefschetz trace
formula.

\bigskip
\bigskip
\noindent {\bf 5. Motivic volume}
\bigskip
\noindent
Here $X$ is again any
algebraic variety of pure dimension $d$.
\bigskip \noindent
{\bf
5.1. Definition.} The {\it motivic volume} of $X$ is
$\mu({\Cal L}(X)) \in \hat {\Cal M}_{\Bbb C}$, thus the motivic
measure of the whole arc space of $X$. Recall that $\mu({\Cal
L}(X)) = \lim_{n \rightarrow \infty} \frac{[\pi_n({\Cal
L}(X))]}{\Bbb L^{nd}}$, and that it equals $[X]$ when $X$ is
nonsingular.
\bigskip \noindent We computed in (3.3) the motivic volume of $X =
\{ xy = 0 \}$ as $\mu({\Cal L}(X)) = 2 \Bbb L$ by the defining
limit procedure. For more complicated $X$, the following formula
in terms of a suitable resolution of singularities is very useful.
\bigskip
\proclaim {5.2. Theorem {\rm[DL3]}} Let $h : Y
\rightarrow X$ be log resolution of $X$; i.e. $h$ is a proper
birational morphism from a nonsingular $Y$ such that the
exceptional locus $Exc$ of $h$ is a normal crossings divisor.
Assume also that the image of $h^\ast(\Omega^d_X)$ in $\Omega^d_Y$
is locally principal, i.e. locally generated by one element.

Denote by $E_i, i \in S$, the irreducible components of $Exc$, and
let $\rho_i - 1$ be the multiplicity along $E_i$ of the divisor
associated to $h^\ast (\Omega^d_X)$, i.e. the (effective) divisor
locally given by the zeroes of a generator of
$h^\ast(\Omega^d_X)$. Finally, set $E^\circ_I := (\cap_{i \in I}
E_i) \setminus (\cup_{\ell \not \in I } E_\ell)$ for $I
\subset S$. Then
$$\mu({\Cal L}(X)) = \sum_{I \subset S} [E^\circ_I] \prod_{i \in
I} \frac{\Bbb L-1}{\Bbb L^{\,\rho_i}-1} = \sum_{I \subset S}
[E^\circ_I] \prod_{i \in I} \frac{1}{[\Bbb P^{\rho_i - 1}]}$$ in
$\hat {\Cal M}_{\Bbb C}$; in particular $\mu({\Cal L}(X))$ belongs
to the subring of $\hat{\Cal M}_{\Bbb C}$, obtained from (the
image of) ${\Cal M}_{\Bbb C}$ by inverting the elements $1 + \Bbb
L + \cdots + \Bbb L^j = [ \Bbb P^j ]$. \endproclaim

We will denote this subring by ${\Cal M}_{loc}$.
\bigskip
{\bf 5.3.} {\sl Example.} Let $X = \{ y^2 - x^3 = 0
\}$ in $\Bbb A^2$. We take $\Bbb A^1 \rightarrow X : u \mapsto
(u^2,u^3)$ as a log resolution. Since $\Omega^1_X$ is generated by
$dx$ and $dy$ (subject to the relation $2y dy = 3x^2 dx$), one
easily verifies that $h^\ast \Omega^1_X$ is generated by $u du$.
Hence the image of $h^\ast \Omega^1_X$ in $\Omega^1_Y$ is
principal and we can apply Theorem 5.2.

Note that $Exc = E_1 = \{ 0 \}$, occurring with multiplicity 1 in
the divisor of $udu$. So $\rho_1 = 2$ and
$$\mu ({\Cal L}(X)) = \Bbb L-1 + \frac{1}{[\Bbb P^1]} =
\frac{\Bbb L^2}{\Bbb L+1}.$$
(Recall that $[X] = \Bbb L$.)

\bigskip
\noindent {\bf 5.4.} {\sl Example.} Let $X = \{ z^2 = xy \}$ in
$\Bbb A^3$.

\noindent {\smc Exercise}. (i) Verify that $\mu({\Cal L}(X)) =
\Bbb L^2$. (The `obvious' log resolution satisfies the assumption
of Theorem 5.2, and the unique component $E_1$ of the exceptional
locus has $\rho_1 = 2$.)

(ii) Note that also $[X] = \Bbb L^2$; this could be interpreted as
the singularity of $X$ being `very mild'.
\bigskip
\noindent {\bf 5.5.} {\smc Exercise}. Compute again that the
motivic volume of $X = \{xy = 0 \}$ is $2 \Bbb L$; now using
Theorem 5.2. (Note here that $[X] = 2 \Bbb L - 1$; one could say
that the motivic volume counts the double point twice.)
\bigskip \noindent
{\bf 5.6.} Recall that for nonsingular $X$ its universal Euler
characteristic $[X] \in K_0 (Var_{\Bbb C})$ specializes to its
Hodge-Deligne polynomial $H(X) \in \Bbb Z[u,v]$ and further to
$\chi(X)  \in \Bbb Z$.

Since $\chi(\cdot)$ and $H(\cdot)$ factor through the image of
${\Cal M}_{\Bbb C}$ in $\hat{\Cal M}_{\Bbb C}$, they induce
natural maps $\chi : {\Cal M}_{loc} \rightarrow \Bbb Q$ and $H :
{\Cal M}_{loc} \rightarrow \Bbb Z[[u,v]]$. Applying these
specialization maps to the motivic measure of $X$ yields new
(numerical) singularity invariants, which generalize the usual
$\chi(X)$ and $H(X)$ for nonsingular $X$. Denef and Loeser call
$\chi(\mu({\Cal L}(X)))$ the {\it arc-Euler characteristic} of $X$.

\medskip
For example the arc-Euler characteristic of $\{ y^2 - x^3 = 0 \}$
is $\frac 1 2$ and the one of $\{ xy = 0 \}$ is $2$.

\bigskip
\bigskip
\noindent {\bf 6. Motivic zeta functions}
\bigskip
\noindent In this section $M$ is a nonsingular irreducible
algebraic variety of dimension $m$, and $f : M \rightarrow \Bbb C$
is a non-constant regular function.
\bigskip
\noindent {\bf 6.1.} For each $n \in \Bbb N$ the morphism $f : M
\rightarrow \Bbb A^1 = \Bbb C$ induces a morphism $f_n : {\Cal
L}_n(M) \rightarrow {\Cal L}_n(\Bbb A^1)$.  A point $\alpha \in
{\Cal L}_n(\Bbb A^1)$ corresponds to an element $\alpha(t) \in
K[t]/(t^{n+1})$ for some field $K \supset \Bbb C$; we denote as usual
the largest $e$ such that $t^e$ divides $\alpha(t)$ by $ord_t
\alpha \in \{ 0,1,\cdots,n,+ \infty \}$. We set
$${\Cal X}_n := \{ \gamma \in {\Cal L}_n(M) \mid ord_t f_n(\gamma) =
n \} \quad \text { for } n \in \Bbb N;$$ it is a locally closed
subvariety of ${\Cal L}_n(M)$.
\bigskip
\noindent {\smc Exercise}. Denote $X := \{ f = 0 \}$. Then $[ \Cal
X_n ] = [ {\Cal L_{n-1}}(X)] - [{\Cal L}_n(X)]$ for $n \geq 1$,
and $[{\Cal X}_0] = [M] - [X]$.

\bigskip
\noindent {\bf Definition.} The {\it motivic zeta function} of $f
: M \rightarrow \Bbb C$ is the formal power series
$$Z(T) := \sum_{n \geq 0} [{\Cal X}_n](\Bbb L^{-m} T)^n$$
in ${\Cal M}_{\Bbb C}[[T]]$.

\bigskip
\noindent {\bf 6.2.} Considering the exercise above, it is not a
surprise that for $X := \{ f = 0 \}$ the series $J(T) = \sum_{n
\geq 0} [{\Cal L}_n(X)]T^n$ and $Z(T)$ are equivalent. Indeed, one
easily verifies that
$$J(T) = \frac{Z(\Bbb L^mT) - \Bbb L^m}{\Bbb L^m T - 1}.$$
\bigskip
\noindent {\bf 6.3.} The definition of $Z(T)$ is inspired by the
$p$-adic Igusa zeta function, associated to a polynomial $f \in
\Bbb Z_p[x_1,\cdots,x_m]$, which is defined as
$$Z_p(s) := \int_{{\Bbb Z}^m_p} |f(x)|^s_p |dx|$$
for $s \in \Bbb C, \Re(s) > 0$, and can be rewritten as

$$\split
Z_p(s) & = \sum_{n \geq 0} volume \{ x \in \Bbb Z^m_p \mid ord_p f(x) = n \} p^{-ns} \\
& = \frac{1}{p^m} \sum_{n \geq 0} \# \{ x \in (\Bbb Z/p^{n+1} \Bbb
Z)^m \mid ord_p f(x) = n \} (p^{-m} p^{-s})^n.
\endsplit
$$

\bigskip
\noindent {\bf 6.4.}  {\smc Exercise}. Write $D$ for the
(effective) divisor of zeros of $f$, i.e. $D$ is ``$\{f = 0 \}$
with multiplicities". Then
$$\int_{{\Cal L}(X)} \Bbb L^{-ord_t D} d \mu = Z (\Bbb L^{-1})$$
in $\hat{\Cal M}_{\Bbb C}$, meaning in particular that the
substitution in the right hand side yields a well-defined element
of $\hat{\Cal M}_{\Bbb C}$.

\bigskip
\noindent {\bf 6.5.} As for the motivic volume, there is an
important (similar) formula for $Z(T)$ in terms of a resolution.

\bigskip
\proclaim{Theorem {\rm [DL2]}} Let $h : Y \rightarrow M$ be an
embedded resolution of $\{ f = 0 \}$; i.e. $h$ is a proper
birational morphism from a nonsingular $Y$ such that $h$ is an
isomorphism on $Y \setminus h^{-1} \{ f = 0 \}$  and $h^{-1} \{ f
= 0 \}$ is a normal crossings divisor. Let $E_i, i \in S$, be the
irreducible components of $h^{-1} \{ f = 0 \}$. For $i \in S$ we
denote by $N_i$ the multiplicity of $E_i$ in the divisor of $f
\circ h$ on $Y$, and by $\nu_i - 1$ the multiplicity of $E_i$ in
the divisor of $h^\ast \omega$, where $\omega$ is a local
generator of $\Omega^m_M$. (Equivalently : $div(f \circ h) =
\sum_{i \in S} N_i E_i$ and $K_{Y|M} = \sum_{i \in S} (\nu_i -
1)E_i.)$ Set finally $E^\circ_I := (\cap_{i \in I} E_i) \setminus
(\cup_{\ell \not \in I} E_\ell)$ for $I \subset S$. Then
$$Z(T) = \sum_{I \subset S} [E^\circ_I] \prod_{i \in I}
\frac{(\Bbb L - 1)T^{N_i}}{\Bbb L^{\nu_i} - T^{N_i}};$$ in
particular $Z(T)$ is rational and belongs more precisely to the
subring of ${\Cal M}_{\Bbb C} [[T]]$ generated by ${\Cal M}_{\Bbb
C}$ and the elements $\frac{T^N}{L^\nu - T^N}$, where $\nu,N \in
\Bbb Z_{> 0}$. \endproclaim
\bigskip
\noindent {\bf 6.6.} {\sl Corollaries.}

\smallskip
(i) In the special case that $X = \{ f = 0 \}$ is a hypersurface
this yields the stated rationality of $J(T)$ in (2.12).

(ii) Let $M = \Bbb A^m$ and $f \in \Bbb Z[x_1, \cdots , x_m]$.
Then by a similar formula of Denef [De2] for the $p$-adic Igusa
zeta functions $Z_p(s)$, Theorem 6.5 yields that $Z(T)$
specializes to the $Z_p(s)$ for all $p$ except  a finite number.
See [DL2] for a precise statement. Similarly $J(T)$ specializes to
$J_p(T)$ for all $p$ except a finite number [DL8, Theorem 6.1].

(iii) For any $f : M \rightarrow \Bbb C$ we now explain how $Z(T)$
specializes to the {\it topological zeta function} of $f$. Using
Theorem 6.5 and the notations there, we evaluate $Z(T)$ at $T =
\Bbb L^{-s}$ for any $s \in \Bbb N$; this yields the well-defined
elements
$$\sum_{I \subset S} [E^\circ_I] \prod_{i \in I} \frac{\Bbb
L-1}{\Bbb L^{\nu_i + sN_i} - 1} = \sum_{I \subset S} [E^\circ_I]
\prod_{i \in I} \frac{1}{[\Bbb P^{\nu_i  + sN_i - 1}]}$$ in (the
image in $\hat {\Cal M}_{\Bbb C}$ of) the localization of ${\Cal
M}_{\Bbb C}$ with respect to the elements $[\Bbb P^j]$. Applying
the Euler characteristic specialization map $\chi(\cdot)$ yields
the rational numbers
$$\sum_{I \subset S} \chi(E^\circ_I) \prod_{i \in I}
\frac{1}{\nu_i + sN_i}$$ for $s \in \Bbb N$. The {\it topological
zeta function} $Z_{top}(s)$ of $f$ is the unique rational function
in one variable $s$ admitting the values above for $s \in \Bbb N$.

Without the specialization argument above it is not at all clear
that $Z_{top}(s)$ does not depend on the chosen resolution $h : Y
\rightarrow M$. In fact $Z_{top}(s)$ was first introduced in [DL1],
in terms of a resolution, and $p$-adic Igusa zeta functions and
the Grothendieck-Lefschetz trace formula were needed to prove
independence of the chosen resolution.

\bigskip
\noindent {\bf 6.7.} We
just mention that there is an important generalization of the
motivic zeta function, working over a relative and equivariant
Grothendieck ring; it specializes by a limit procedure to objects
in (an equivariant version of) ${\Cal M}_{\Bbb C}$, which are
shown to be a good {\it virtual motivic incarnation} of the Milnor
fibres of $f$ at the points of $\{ f = 0 \}$. It is quite
remarkable that a definitely non-algebraic notion as the Milnor
fibre has such an algebraic incarnation. See [DL2][DL7].

Moreover these objects satisfy a {\it motivic Thom-Sebastiani
Theorem}, generalizing the known results of Varchenko and Saito.
See [DL4].
\bigskip
\noindent {\bf 6.8.} {\sl Monodromy Conjecture.}
\smallskip
\noindent There is an intriguing conjectural relation between the
poles of the topological zeta function and the eigenvalues of the
local monodromy of $f$.

\medskip
\noindent \proclaim{Monodromy conjecture} If $s_0$ is a pole of
$Z_{top}(s)$, then $e^{2 \pi is_0}$ is an eigenvalue of the local
monodromy action on the cohomology of the Milnor fibre of $f$ at
some point of $\{ f = 0 \}$.
\endproclaim

\medskip
One can also state the analogous conjecture for the motivic zeta
function, but then one has to be careful with the notion of pole,
see [RV2]. Alternatively, we can formulate this monodromy
conjecture for $Z(T)$ as follows, without mentioning poles [DL2] :

\noindent {\sl $Z(T)$ belongs to the ring generated by $\Cal
M_{\Bbb C}$ and the elements $\frac{T^N}{L^\nu - T^N}$, where
$\nu,N \in \Bbb Z_{> 0}$ and $e^{2 \pi i \frac{\nu}N}$ is an
eigenvalue of the local monodromy as above.}

\medskip
\noindent Actually, it was originally stated for the $p$-adic
Igusa zeta function, being even more remarkable, for then it
relates number theoretical invariants of $f \in \Bbb
Z[x_1,\cdots,x_m]$ to differential topological invariants of $f$,
considered as function $\Bbb C^n \rightarrow \Bbb C$.

The conjecture was shown by Loeser for $M = \Bbb A^2$ [Loe1]; a
shorter proof in dimension 2 is in [Ro]. In dimension 3 there is a
lot of `experimental evidence' [Ve1], and by now various special
cases are proved [ACLM1][ACLM2][Loe2][RV1].

\bigskip
\noindent {\sl Example.} Let $M = \Bbb A^2$ and $f = y^2 - x^3$.

\noindent {\smc Exercise}. Compute, using Theorem 6.5,
$$Z(T) = \Bbb L^2 (\Bbb L - 1) \frac{\Bbb L^5 - \Bbb L^3 T + \Bbb
L^3 T^2 - T^5}{(\Bbb L^5 - T^6)(\Bbb L - T)}$$
and
$$Z_{top}(s) = \frac{5 + 4s}{(5 + 6s)(1 + s)}.$$
(This is how we computed $J(T)$ in Example 2.14.) In particular,
the poles of $Z_{top}(s)$ are $-1$ and $-5/6$. On the other hand,
it is well known that the monodromy eigenvalues of $f$ are $1,
e^{\frac{\pi i}{3}}$, and $e^{- \frac{\pi i}{3}}$. Hence the
monodromy conjecture is indeed satisfied here.

\bigskip
\noindent {\sl Note.} The previous example was too simple to exhibit the
`typical' situation. Each irreducible component $E_i$ in Theorem
6.5 induces a candidate-pole $-\frac{\nu_i}{N_i}$, and quite
miraculously, for a generic example with a lot of components
$E_i$, `most' of these candidates cancel.  This experimental fact
is compatible with the monodromy conjecture, see [Ve1].

\bigskip
\bigskip
\noindent {\bf 7. Batyrev's stringy invariants}
\bigskip
\noindent Using motivic integration, Batyrev [Ba1][Ba2] introduced new
singularity invariants for algebraic varieties with `mild'
singularities, more precisely with at worst log terminal
singularities. He used them for instance to formulate a
topological mirror symmetry test for singular Calabi-Yau
varieties, to give a conjectural definition for stringy Hodge
numbers, and to prove a version of the McKay correspondence.

We first explain log terminal and related singularities; for this
we need the Gorenstein notion.

\bigskip
\noindent {\bf 7.1.} Let $X$ be a {\it normal} algebraic variety
of dimension $d$. In particular $X$ is irreducible, $X_{\sing}$ has codimension at least $2$ in $X$, and $X$ has a well
defined canonical divisor $K_X$ (up to linear equivalence). One
can view (a representative of) $K_X$ as the divisor of zeroes and
poles of a {\it rational} differential $d$-form on $X$; it is also
the Zariski-closure of the usual canonical divisor on $X_{\reg}$.

When $X$ is nonsingular, $K_X$ is a Cartier divisor, i.e. locally
given by one equation. This is not true in general.

\bigskip \noindent
{\bf Definition.} A normal variety $X$ is {\it Gorenstein} if
$K_X$ is a Cartier divisor. Alternatively : $X$ is Gorenstein if
the rational differential $d$-forms on $X$, which are regular on
$X_{\reg}$, are locally generated by one element.

\bigskip\noindent
{\sl Example.} Let $X = \{z^2 = xy \}$; then those differential 2-forms
are generated by $\frac{dx \wedge dy}{2z} = \frac{dx \wedge dz}{x} = -
\frac{dy \wedge dz}{y}$ (which is indeed regular on $X_{\reg}$).
\bigskip
\noindent This notion is quite general; for instance all (normal)
hypersurfaces and even complete intersections are Gorenstein.

\bigskip
\noindent {\bf 7.2.} We now introduce a certain `badness' for
singularities, in terms of numerical invariants of a resolution.

Let $X$ be Gorenstein of dimension $d$. Take a log resolution $\pi
: Y \rightarrow X$ of $X$ and denote by $E_i, i \in S$, the
irreducible components of the exceptional locus $Exc$ of $h$. We
associate as follows an integer $a_i$ to each $E_i$.
\medskip
(1) {\it Description with divisors.} Since $K_X$ is Cartier, the
pullback $\pi^\ast K_X$ makes sense and one can consider the
relative canonical divisor $K_{Y|X} = K_Y - \pi^\ast K_X$, which
is supported on $Exc$. Then $a_i - 1$ is the multiplicity of $E_i$
in $K_{Y|X}$, i.e. $K_{Y|X} = \sum_{i \in S} (a_i - 1)E_i$.
\medskip
(2) {\it Description with differential forms.} Take a general
point $Q_i$ of $E_i$ and local coordinates $y_1,y_2,\cdots,y_d$
around $Q_i$ such that the local equation of $E_i$ is $y_1 = 0$.
Let $\omega_i$ be a local generator around $\pi(Q_i)$ of the
$d$-forms on $X$, which are regular on $X_{\reg}$. (Such an
$\omega_i$ exists by the Gorenstein property.) Then around $Q_i$
one can write $\pi^\ast \omega_i$ as
$$\pi^\ast \omega_i = u y^{a_i-1}_1 dy_1 \wedge dy_2 \wedge \cdots
\wedge dy_d ,$$ where $u$ is regular and nonzero around $Q_i$.

\medskip
\noindent In general the $a_i \in \Bbb Z$, and when $X$ is
nonsingular they satisfy $a_i \geq 2$.
\medskip \noindent
{\sl Terminology.} One calls $a_i$ the {\it log discrepancy} of $E_i$
with respect to $X$ (and $a_i - 1$ the {\it discrepancy}).

\bigskip
\noindent {\sl Example.} The standard log resolution of
$X=\{z^2=xy\}$ has one exceptional curve $E \cong \Bbb P^ 1$ with
log discrepancy $a=1$.

\bigskip
\noindent {\bf 7.3.} We also have to consider a technical
generalization: a normal variety is called {\it $\Bbb
Q$-Gorenstein} if $rK_X$ is Cartier for some $r \in \Bbb Z_{>
0}$. Then the log discrepancies are defined analogously by
$K_{Y|X} = \sum_{i \in S} (a_i - 1)E_i$, which should be
considered as an abbreviation of $rK_{Y|X} = rK_Y - rK_X = \sum_{i
\in S} r(a_i - 1)E_i$. Now the $r(a_i - 1) \in \Bbb Z$, and hence
$a_i \in \frac 1 r \Bbb Z$.

\bigskip
\noindent {\sl Example.} Let $X$ be the quotient of $\Bbb A^2$ by
the action of $\mu_3=\{z\in \Bbb C \mid z^3=1\}$ given by
$(x,y)\mapsto(\epsilon x,\epsilon y)$ for $\epsilon \in \mu_3$.
Concretely, $X$ is given in $\Bbb A^4$ by the equations
$$
\{u_1u_3-u_2^2=u_2u_4-u_3^2=u_1u_4-u_2u_3=0\}\, ,
$$
in particular it is {\sl not} a complete intersection. Here $K_X$
is not Cartier; a representative of $K_X$ is for example
$\{u_1=u_2=u_3=0\}$. However, $3K_X$ is Cartier; a representative
is $\{u_1=0\}$.

The standard log resolution of $X$ has one exceptional curve
$E\cong \Bbb P^1$ with log discrepancy $a=\frac 23$.

\medskip
A nice introduction to these notions is in [Re1].
\bigskip
\noindent {\bf 7.4.} {\bf Definition.} (i) Let $X$ be a $\Bbb
Q$-Gorenstein variety. Take a log resolution $\pi : Y \rightarrow
X$ of $X$; let $E_i, i \in S$, be the irreducible components of
the exceptional locus of $\pi$ with log discrepancies $a_i$. Then
$X$ is called {\it terminal, canonical, log terminal} and {\it log
canonical} if $a_i > 1, a_i \geq 1, a_i > 0$ and $a_i \geq 0$,
respectively, for all $i \in S$.
\smallskip \noindent One can show
that these conditions do not depend on the chosen resolution.

\smallskip
(ii) We say that $X$ is {\it strictly log canonical} if it is log
canonical but not log terminal.
\bigskip
\noindent We should note  that $0$ is indeed the relevant `border
value' here; if some $a_i < 0$ on some log resolution, then one
can easily construct log resolutions with arbitrarily negative
$a_i$.

The log terminal singularities should be considered `mild', the
singularities which are not log canonical `general', and the
strictly log canonical ones as a special `border' class.
\bigskip
\noindent {\bf 7.5.} {\sl Example.}
 (1) When $X$ is a surface ($d=2$) terminal is equivalent to
 non-singular, the canonical singularities are precisely the
 so-called ADE singularities or rational double points, and the
 log terminal singularities are precisely the Hirzebruch-Jung or quotient
 singularities.

\smallskip
(2) Let $X = \{ x^k_1 + x^k_2 + \cdots + x^k_{d+1} = 0 \}$ in
$\Bbb A^{d+1}$. The origin is the only singular point of $X$, and
the blowing-up with the origin as centre yields a log resolution
$\pi : Y \rightarrow X$ of $X$ with exceptional locus consisting
of one irreducible component $E$, which is isomorphic to $\{x^k_1
+ x^k_2 + \cdots + x^k_{d+1} = 0 \} \subset \Bbb P^d$.
\smallskip
\noindent {\smc Exercise}. (i) The log discrepancy of $E$ with
respect to $X$ is $d + 1 - k$.

(ii) $X$ is log terminal, strictly log canonical, and not log
canonical when $k < d + 1, k = d + 1$, and $k > d + 1$,
respectively.
\bigskip
\noindent {\bf  7.6.}  There are nice results of Ein, Musta\c t\v a
and Yasuda, relating the previous notions with jet spaces.
\bigskip
\noindent \proclaim{Theorem {\rm [Mu1][EMY][EM]}} Let $X$ be a
normal variety, which is locally a complete intersection. Then $X$
is terminal, canonical, and log canonical if and only if ${\Cal
L}_n(X)$ is normal, irreducible, and equidimensional,
respectively, for every $n$. \endproclaim
\bigskip
\noindent {\bf 7.7.} {\bf Definition.} Let $X$ be a log terminal
algebraic variety. Take a log resolution $\pi : Y \rightarrow X$
of $X$. Let $E_i, i \in S$, be the irreducible components of the
exceptional locus of $\pi$ with log discrepancies $a_i$ ($\in \Bbb Q_{>0}$). Denote
also $E^\circ_I := (\cap_{i \in I} E_i) \setminus (\cup_{\ell \not
\in I} E_\ell)$ for $I \subset S$.

 (i) The {\it stringy Euler number} of $X$ is
 $$e_{st}(X) := \sum_{I \subset S} \chi(E^\circ_I) \prod_{i \in I}
 \frac{1}{a_i}.$$

  (ii) The {\it stringy E-function} of $X$ is
 $$E_{st} (X) := \sum_{I \subset S} H(E^\circ_I) \prod_{i \in I}
 \frac{uv-1}{(uv)^{a_i}-1}.$$

 (iii) The {\it stringy ${\Cal E}$-invariant} of $X$ is
 $${\Cal E}_{st} (X) := \sum_{I \subset S} [E^\circ_I] \prod_{i \in I}
 \frac{\Bbb L - 1}{\Bbb L^{a_i} - 1}.$$
\medskip
\noindent {\sl Remarks.} (1) Clearly $e_{st} (X) \in \Bbb Q;
E_{st}(X)$ is a rational function in $u,v$ (with `fractional
powers'), and ${\Cal E}_{st}(X)$ lives in a finite extension of
$\hat{\Cal M}_{\Bbb C}$. We have specialization maps ${\Cal
E}_{st}(X) \mapsto E_{st} (X) \mapsto e_{st} (X)$.

(2) Strictly speaking, Batyrev defined and used only the levels (i) and (ii)
[Ba2][Ba3].
\bigskip \noindent
When $X$ is nonsingular, ${\Cal E}_{st}
(X) = [X]$ (this is 4.2), and of course $E_{st} (X) = H(X)$ and
$e_{st}(X) = \chi(X)$. So also these invariants are new
singularity invariants, generalizing $[ \cdot ], H(\cdot)$ and
$\chi(\cdot)$, respectively, for nonsingular $X$. (Just as the
motivic volume and its specializations. We give a comparing
example in 7.11.)
\bigskip
\noindent {\bf 7.8.} The crucial point is that the defining
expressions above do not depend on the chosen resolution. We
indicate three different arguments, supposing for simplicity that
$X$ is Gorenstein, i.e. the $a_i \in \Bbb Z_{> 0}$.
\bigskip
\noindent (1) Let $\pi : Y \rightarrow X$ and $\pi^\prime :
Y^\prime \rightarrow X$ be two log resolutions of $X$. By the
formula of Proposition 3.6 we have in fact
$$\sum_{I \subset S} [E^\circ_I] \prod_{i \in I} \frac{\Bbb
L-1}{\Bbb L^{a_i}-1} = \int_{{\Cal L}(Y)} \Bbb L^{-ord_t K_{Y|X}}
d \mu.$$
So we must show that $\int_{{\Cal L}(Y)} \Bbb L^{-ord_t
K_{Y|X}} d \mu = \int_{{\Cal L}(Y^\prime)} \Bbb L^{-ord_t
K_{Y^\prime|X}} d \mu$. To this end we take a log resolution $\rho
: Z \rightarrow X$, dominating $\pi$ and $\pi^\prime$; i.e. we
have $\rho : Z \overset \sigma \to \rightarrow Y \overset \pi \to
\rightarrow X$ and $\rho : Z \overset \sigma^\prime \to
\rightarrow Y^\prime \overset \pi^\prime \to \rightarrow X$. By
the change of variables formula in (3.8) we have
$$\int_{{\Cal
L}(Y)} \Bbb L^{-ord_t K_{Y|X}} d \mu = \int_{{\Cal L}(Z)} \Bbb
L^{-ord_t(\sigma^\ast K_{Y|X} + K_{Z|Y})} d \mu = \int_{{\Cal
L}(Z)} \Bbb L^{-ord_t(K_{Z|X})} d \mu \, ,$$
 and of course the same is
true for the integral over ${\Cal L}(Y^\prime)$.

This is essentially Batyrev's proof.
\bigskip
\noindent (2) We can define ${\Cal E}_{st}(X)$ intrinsically,
using motivic integration on $X$ [Ya][DL6]. There is an ideal sheaf
${\Cal I}_X$ on $X$ such that
$${\Cal E}_{st}(X) = \int_{{\Cal L}(X)} \Bbb L^{ord_t {\Cal I}_X
} d \mu,$$ using the setting of (3.5) and (3.7). More precisely,
denoting by $\omega_X$ the sheaf of differential $d$-forms on $X$
which are regular on $X_{\text{reg}}$, we have a natural map
$\Omega^d_X \rightarrow \omega_X$ whose image is ${\Cal I}_X
\omega_X$. See [Ya, Lemma 1.16].
\bigskip
\noindent (3) Using the Weak Factorization Theorem, see below, one
essentially has to show that the defining expressions in (7.7) do
not change after blowing-up $Y$ in a nonsingular centre which
intersects $\cup_{i \in S} E_i$ transversely. This is
straightforward.
\medskip
\proclaim {7.9. Weak Factorization Theorem {\rm [AKMW][W\l]}}

(1) Let $\phi : Y -\! \rightarrow Y^\prime$ be a proper birational map between nonsingular irreducible varieties, and let $U \subset Y$ be an open set where $\phi$ is an isomorphism.  Then $\phi$ can be factored as follows into a sequence of blow--ups and blow--downs with smooth centres disjoint from $U$.

There exist nonsingular irreducible varieties $Y_1, \dots , Y_{\ell -1}$ and a sequence of birational maps
$$Y = Y_0 - \overset \phi_1 \to \rightarrow Y_1 - \overset \phi_2 \to \rightarrow \cdots -\! \overset \phi_{i-1} \to \rightarrow Y_{i-1} - \overset \phi_i \to \rightarrow Y_i - \! \overset \phi_{i+1} \to \rightarrow \cdots - \! \overset \phi_{\ell - 1} \to \rightarrow Y_{\ell - 1} -  \overset \phi_\ell \to \rightarrow Y_\ell = Y^\prime$$
where $\phi = \phi_\ell \circ \phi_{\ell - 1} \circ \cdots \circ \phi_2 \circ \phi_1$, such that each $\phi_i$ is an isomorphism over $U$ (we identify $U$ with an open in the $Y_i$), and for $i = 1, \dots , \ell$ either $\phi_i : Y_{i-1} -\! \rightarrow Y_i$ or $\phi^{-1}_i : Y_i -\! \rightarrow Y_{i-1}$ is the blowing--up at a nonsingular centre disjoint from $U$, and is thus a morphism.

(1$^\prime$) There is an index $i_0$ such that for all $i \leq i_0$ the map $Y_i \rightarrow Y$ is a morphism, and for $i \geq i_0$ the map $Y_i \rightarrow Y^\prime$ is a morphism.

(2) If $Y \setminus U$ and $Y^\prime \setminus U$ are normal crossings divisors, then the factorization above can be chosen such that the inverse images of these divisors under $Y_i \rightarrow Y$ or $Y_i \rightarrow Y^\prime$ are also normal crossings divisors, and such that the centres of blowing--up of the $\phi_i$ or $\phi^{-1}_i$ intersect these divisors transversely.

\endproclaim
\medskip
\noindent
{\sl Remark.}  (i) In [AKMW] and [W{\l}] the theorem is stated for a birational map $\phi$ between complete $Y$ and $Y^\prime$; the generalization to proper birational maps between not necessarily complete $Y$ and $Y^\prime$ is mentioned by Bonavero [Bo].

(ii) In [AKMW, Theorem 0.3.1] the first claim of (2) is not explicitly stated, but can be read off from the proof (see [AKMW, 5.9 and 5.10]).


\bigskip
\noindent {\bf 7.10.} {\sl Important Intermezzo.} Using weak
factorization instead of motivic integration, we can define ${\Cal
E}_{st}(X)$ in a localization of (a finite extension of) ${\Cal
M}_{\Bbb  C}$, which is a priori finer than in (a finite extension
of) $\hat {\Cal M}_{\Bbb C}$, since we do not know whether the
natural map $\Cal M_{\Bbb C} \rightarrow \hat {\Cal M}_{\Bbb C}$
is injective.
\medskip
This remark also applies e.g. to (4.1), yielding $[X] = [Y]$ in
the localization of ${\Cal M}_{\Bbb C}$ with respect to the $[\Bbb
P^j]$ instead of merely in $\hat{\Cal M}_{\Bbb C}$.
\bigskip
\noindent {\bf 7.11.} {\sl Example.} Let $X = \{ x^k_1 + x^k_2 +
\cdots + x^k_{d+1} = 0 \} \subset \Bbb A^{d+1}.$

\noindent {\smc Exercise}. We use the notation $E$ of Example 7.5.

(i) ${\Cal E}_{st} (X) = (\Bbb L - 1)[E] + [E] \frac{\Bbb
L-1}{\Bbb L^{d+1-k}-1} ,$

(ii) $\mu({\Cal L}(X)) = (\Bbb L - 1)[E] + [E] \frac {\Bbb
L-1}{\Bbb L^d - 1} ,$

(iii) $[X] = (\Bbb L-1)[E] + 1.$
\medskip
\noindent (Note also that (ii) and (iii) are consistent with
Example 5.4.)
\bigskip
\noindent {\bf 7.12.} {\sl Applications.}

(i) Topological mirror symmetry test for singular Calabi-Yau
mirror pairs [Ba2].

(ii) A conjectural definition of stringy Hodge numbers for certain
canonical Gorenstein varieties [Ba2].

(iii) A proof of a version of the McKay correspondence
[Ba3][DL6][Ya1].

(iv) A new birational invariant for varieties of nonnegative
Kodaira dimension, assuming the Minimal Model Program [Ve2, (2.8)].

\bigskip
\bigskip
\noindent {\bf 8. Stringy invariants for general singularities}
\bigskip
\noindent In this section $X$ is a $\Bbb Q$-Gorenstein  variety.
\bigskip
\noindent {\bf 8.1.} For a log resolution $\pi : Y \rightarrow X$ of
$X$, we use the notations $E_i$ and $a_i, i \in S$, and
$E^\circ_I, I \subset S$, as before. There are (at least) two
natural questions concerning a possible generalization of
Batyrev's stringy invariants beyond the log terminal case.
\bigskip
\noindent {\smc Question I.} Suppose there exists at least one log
resolution $\pi : Y \rightarrow X$ of $X$ for which {\it all} log
discrepancies $a_i \ne 0$. Is (e.g.)
$$\sum_{I \subset S} \chi (E^\circ_I) \prod_{i \in I}
\frac{1}{a_i}$$ independent of a chosen such resolution ?

\medskip
This question is still open (a positive answer would yield a
generalized stringy invariant for those $X$ admitting such a log
resolution). Note that, when using the weak factorization theorem
to connect two such log resolutions by chains of blowing-ups, log
discrepancies on \lq intermediate varieties\rq\ could be zero,
obstructing an obvious attempt of proof.
\bigskip
\noindent {\smc Question II.} Do there exist any kind of
invariants, associated to all or `most' $\Bbb Q$-Gorenstein
varieties, which coincide with Batyrev's stringy invariants if the
variety is log terminal ?
\medskip
Concerning this question, we obtained the following result [Ve4].
We associated invariants to `almost all' $\Bbb Q$-Gorenstein
varieties, more precisely to all $\Bbb Q$-Gorenstein varieties
without strictly log canonical singularities, which do generalize
Batyrev's invariants for log terminal varieties. (Note that in
particular log discrepancies can be zero in a log resolution of a
non log canonical variety !)

\item{$\bullet$} To construct these invariants we have to assume
Mori's Minimal Model Program (in fact the relative and log
version). \item{$\bullet$} As in the previous section, we can work
on any level : $\chi(\cdot), H(\cdot)$, and $[\cdot]$.  For
simplicity we treat here just the roughest level $\chi(\cdot)$;
the other levels are analogous.
\bigskip
\noindent {\bf 8.2.} We associate to any $\Bbb Q$-Gorenstein $X$ without
strictly log canonical singularities a rational function $z_{st}
(X;s)$ in one variable $s$, the {\it stringy zeta function}  of
$X$. It will turn out that for log terminal $X$, this rational
function is in fact a constant and equal to $e_{st}(X)$.
\bigskip
 We just present the main idea of our construction. The
`pragmatic' idea is to split the log discrepancies $a_i$ of a log
resolution $\pi : Y \rightarrow X$ as $a_i = \nu_i + N_i$ such
that $(\nu_i,N_i) \ne (0,0)$ for all $i$, and to define
$z_{st}(X;s)$ as
$$\sum_{I \subset S} \chi (E^\circ_I) \prod_{i \in I}
\frac{1}{\nu_i + sN_i} \in \Bbb Q(s).$$ This is done in a
geometrically meaningful way via factoring $\pi$ through a certain
`partial resolution' $p : X^m \rightarrow X$ of $X$, which is
called a {\it relative log minimal model of} $X$. This is a
natural object in the (relative, log) Minimal Model Program;
important here is that it is not unique and that $X^m$ can have
certain mild singularities.

For the specialists : $p$ is a proper birational morphism, $X^m$
is $\Bbb Q$-factorial, the pair $(X^m, E^m)$ is divisorial log
terminal, and $K_{X^m} + E^m$ is $p$-nef, where $E^m$ denotes the
reduced exceptional divisor of $p$. References for these notions are e.g. in [KM][KMM][Ma].
\bigskip
We consider the factorization $\pi : Y \overset h \to \rightarrow
X^m \overset p \to \rightarrow X$. In general $h$ is only a
birational map (maybe not everywhere defined), but we suppose for
the moment that it is a morphism. We justify this later. Denoting
as usual by $E_i, i \in S$, the irreducible components of the
exceptional divisor of $\pi$, we let $E^m_i, i \in S^m$, be the
images in $X^m$ of those $E_i$ which `survive' in $X^m$, i.e.
which are not contracted by $h$ to varieties of smaller dimension.
Then
$$\split
\sum_{i \in S} a_i E_i & = K_Y + \sum_{i \in S} E_i - \pi^\ast K_X
\\
& = \undersetbrace (1) \to{K_Y + \sum_{i \in S} E_i - h^\ast
(K_{X^m} + \sum_{i \in S^m} E^m_i)} + \undersetbrace (2)
\to{h^\ast(K_{X^m} + \sum_{i \in S^m} E^m_i) - h^\ast p^\ast K_X}.
\endsplit
$$
Both (1) and (2) are divisors on $Y$, supported on $\cup_{i \in S}
E_i$. We write (1) as $\sum_{i \in S} \nu_i E_i$; all $\nu_i \geq
0$ because the pair $(X^m,\sum_{i \in S^m} E^m_i)$ has only mild
singularities (more precisely, because it is divisorial log
terminal). We can rewrite (2) as
$$h^\ast(K_{X^m} + \sum_{i \in S^m} E^m_i - p^\ast K_X) = h^\ast
(\sum_{i \in S^m} a_i E^m_i);$$ and it is well known that all
$a_i, i \in S^m$, are non-positive (more precisely, this follows
since $K_{X^m} + \sum_{i \in S^m} E^m_i$ is $p$-nef). So we can
write (2) as $\sum_{i \in S} N_i E_i$ where all $N_i \leq 0$.
\medskip
With these definitions of $\nu_i$ and $N_i$ we indeed have $a_i =
\nu_i + N_i$ for $i \in S$, with moreover $\nu_i \geq 0$ and $N_i
\leq 0$. One can show that, if $X$ has no strictly log canonical
singularities, the situation $\nu_i = N_i = 0$ {\sl cannot} occur.

When $X$ is log terminal, the morphism $p : X^m \rightarrow X$ has
{\sl no} exceptional divisors, so $S^m = \emptyset$, all $N_i = 0$
and $\nu_i = a_i$, and as promised $z_{st} (X;s) = e_{st}(X)$.

\bigskip
In fact we FIRST choose a relative log minimal model $p :
X^m \rightarrow X$ of $X$, we secondly choose a log resolution $h
: Y \rightarrow X^m$ of the pair $(X^m,E^m)$, where $E^m$ is the
reduced exceptional divisor of $p$, and then we put $\pi := p \circ h$.

The point is again that $z_{st}(X;s)$ is independent of both
choices, for which a crucial ingredient is the Weak Factorization
Theorem.
\bigskip
\proclaim {8.3. Theorem {\rm[Ve4]}} Let $X$ be any {\rm
surface} without strictly log canonical singularities. Then
$$\lim_{s \rightarrow 1} z_{st} (X;s) \in \Bbb Q.$$
\endproclaim
\medskip
\noindent (Recall that this is non-obvious since some $a_i$ can be
zero. The clue is that if $a_i=0$, then $E_i$ must be rational and
must intersect exactly once or twice other components; this then
easily implies the cancellation of $\nu_i + sN_i$ in the
denominator of $z_{st}(X;s)$.) So we can define in dimension 2 a
generalized stringy Euler {\sl number} $e_{st}(X)$ as the limit
above for any such surface $X$. In fact we constructed this
generalized $e_{st}(X)$ in [Ve3] by a `direct' approach.


\bigskip
\centerline{\beginpicture
 \setcoordinatesystem units <.5truecm,.5truecm>
 \put {$(i)$} at 13 0
 \put {$(1)$} at 13 3
 \put {$(k)$} at 13 -3
 \plot 2 3  -2 0  2 -3 /
 \putrule from -2 0 to 5 0
 \putrule from 7 0 to 10 0
 \putrule from 2 3 to 5 3
 \putrule from 7 3 to 10 3
 \putrule from 2 -3 to 5 -3
 \putrule from 7 -3 to 10 -3
 \multiput {\dots} at 6 0  6 3  6 -3    /
 \multiput {\vdots} at .5 1.2  .5 -.8  /
 \multiput {$\bullet$} at  -2 0  2 0  4 0  8 0  10 0  2 3  4 3  8 3  10 3  2 -3  4 -3  8 -3  10 -3 /
 \put {$E$} at -2.6 0
 \put {$E_{r_i}^{(i)}$} at 2 -.8
 \put {$E_{r_{i}-1}^{(i)}$} at 4 -.8
 \put {$E_2^{(i)}$} at 8 -.8
 \put {$E_1^{(i)}$} at 10 -.8
\endpicture}
\vskip .7truecm \centerline{\smc Figure 1} \vskip 1truecm
\noindent {\bf 8.4.} {\sl Example} [Ve3]. Let $P \in X$ be a
normal surface singularity with dual graph of its minimal log
resolution $\pi : X \rightarrow S$ as in Figure 1.
 There is a central curve $E$ with genus $g$ and
self-intersection number $-\kappa$, and all other curves are
rational. Each attached chain $E^{(i)}_1 - \dots - E^{(i)}_{r_i}$
is determined by two co-prime numbers $n_i$ and $q_i$, which are
the absolute value of the determinant of the intersection matrix
of $E^{(i)}_1, \dots , E^{(i)}_{r_i}$ and $E^{(i)}_1, \dots ,
E^{(i)}_{r_i-1}$, respectively. Finally, we denote by $d$ the
absolute value of the determinant of the total intersection matrix
of $\pi^{-1}P$.
 This is a quite large class of
singularities; it includes all weighted homogeneous isolated
complete intersection singularities, for which the numbers $\{ g;
\kappa; (n_1,q_1), \cdots , (n_k,q_k) \}$ are called the {\it
Seifert invariants} of the singularity.

\medskip If $P \in X$ is not strictly log canonical, then

$$e_{st}(X) = \lim_{s\to 1} z_{st}(X;s) = \frac 1a (2 - 2g - k + \sum^k_{i=1} n_i) + \chi(X\setminus \{P\})\, ,
$$
where
$$a = \frac{2-2g-k + \sum^k_{i=1} \frac{1}{n_i}}{\kappa - \sum^k_{i=1} \frac{q_i}{n_i}}
= \frac{\prod^k_{i=1} n_i}{d} (2 - 2g - k + \sum^k_{i=1}
\frac{1}{n_i})$$
 is the log discrepancy of $E$.

We note that some other log
 discrepancies might be zero.
 A particular example is the so-called {\sl triangle singularity},
 given by $g=0, \kappa=1, k=3$ and $r_1=r_2=r_3=1$. So, concretely,
 there is a central rational curve with self-intersection $-1$ to
 which three other rational curves are attached. Then $a=-1$ and
 the three other log discrepancies are zero, and
$e_{st}(X) = 1- (n_1+n_2+n_3) + \chi(X\setminus \{P\})$.

 \medskip
When such $P \in X$ is a weighted homogeneous isolated {\sl
hypersurface} singularity, this generalized stringy Euler number
appears in some Taylor expansion associated to it, studied by
N\'emethi and Nicolaescu [NN].

\bigskip\noindent
{\bf 8.5.} {\sl Example.} [Ve4] Here we mention a concrete example
of a threefold singularity $P \in X$, having an exceptional
surface with log discrepancy zero in a log resolution, and such
that nevertheless $\lim_{s\to 1} z_{st}(X;s) \in \Bbb Q$, i.e.
such that the evaluation $z_{st}(X;1)$ makes sense.

\smallskip
Let $X$ be the hypersurface $\{x^4+y^4+z^4+t^5=0\}$ in $\Bbb A^4$;
its only singular point is $P=(0,0,0,0)$. We sketch the following
constructions in Figure 2; we denote varieties and their strict
transforms by the same symbol.

The blowing-up $\pi_1:Y_1\to X$ with centre $P$ is already a
resolution of $X$ ($Y_1$ is smooth). Its exceptional surface $E_1$
is the affine cone over the smooth projective plane curve
$C=\{x^4+y^4+z^4=0\}$. Let $\pi_2:Y_2\to Y_1$ be the blowing-up
with centre the vertex $Q$ of this cone, and exceptional surface
$E_2 \cong \Bbb P^2$. Then $E_1 \subset Y_2$ is a ruled surface
over $C$ which intersects $E_2$ in a curve isomorphic to $C$. The
composition $\pi=\pi_1 \circ \pi_2$ is a log resolution of $P\in
X$, and one easily verifies that the log discrepancies are $a_1=0$
and $a_2=-1$; in particular $P\in X$ is not log canonical.

Now $E_1 \subset Y_2$ can be contracted (more precisely one can
check that the numerical equivalence class of the fibre of the
ruled surface $E_1$ is an extremal ray). Let $h:Y_2\to X^m$ denote
this contraction, and let $\pi=p\circ h$. As the notation
suggests, one can verify that $K_{X^m}+E_2$ is $p$-nef, implying
that $(X^m, E_2)$ is a relative log minimal model of $P\in X$.

\vskip 1true cm \centerline{\beginpicture \setcoordinatesystem
units <.45truecm,.45truecm> \putrectangle corners at 6 14 and 16 4
\multiput {$E_2$} at 15 7.7  25 2 / \multiput {$E_1$} at 8.3 13
1.6 4.2  / \multiput {$C$} at 11.5 8.8   23.5 -.2  /
\ellipticalarc axes ratio 3:1.1  360 degrees from 13 12 center at
11 12 \ellipticalarc axes ratio 3:1.1  360 degrees from 13 9
center at 11 9 \ellipticalarc axes ratio 3:1.1  360 degrees from
15 9 center at 11 9 \ellipticalarc axes ratio 3:1.1  360 degrees
from 13 6 center at 11 6 \putrule from 13 6 to 13 12 \putrule from
9 6 to 9 12 \putrectangle corners at -3 -5 and 3 5 \ellipticalarc
axes ratio 3:1.1  360 degrees from 2 3 center at 0 3
\ellipticalarc axes ratio 3:1.1  360 degrees from 2 -3 center at 0
-3 \setquadratic \plot  2 3  1.5 1.5  0 0  -1.5 -1.5  -2 -3 /
        \plot  -2 3  -1.5 1.5  0 0  1.5 -1.5  2 -3 /
\put {$\bullet$} at 0 0 \put {$Q$} at 1 0 \putrectangle corners at
18 -3 and 28 3 \ellipticalarc axes ratio 3:1.1  360 degrees from
25 0 center at 23 0 \ellipticalarc axes ratio 3:1.1  360 degrees
from 27 0 center at 23 0 \putrectangle corners at 8 -3 and 14 -7
\put {$\bullet$} at 11 -5 \put {$P$} at 11.8 -5 \arrow <.3truecm>
[.2,.6] from 5 8 to 2 6 \arrow <.3truecm> [.2,.6]  from 17.5 6 to
20.5 4 \arrow <.3truecm> [.2,.6]  from 11 2.5 to 11 -1.5 \arrow
<.3truecm> [.2,.6] from 4 -2 to 7 -5 \arrow <.3truecm> [.2,.6]
from 17 -2 to 15 -4 \put {$\pi_2$} at  3.5 7.9 \put {$h$} at  19.1
5.8 \put {$\pi$} at  11.7 .4 \put {$\pi_1$} at  5.7 -2.4 \put
{$p$} at  16 -2 \put {$Y_2$} at  16.9 13 \put {$Y_1$} at -3.7 0
\put {$X$} at 14.7 -6 \put {$X^m$} at 29 0
\endpicture}
\vskip 1truecm \centerline{\smc Figure 2} \vskip 1truecm

\noindent Denoting as usual
$$
K_{Y_2}=h^*(K_{X^m}+E_2)+(\nu_1-1)E_1+(\nu_2-1)E_2
\quad\text{and}\quad h^*(a_2E_2)= N_1E_1+N_2E_2
$$
we have clearly that $\nu_2=0$ and $N_2=-1$, and one computes that
$\nu_1=\frac 15$ and $N_1=-\frac 15$. So
$$\align
z_{st}(X;s)&=\frac{\chi(C)}{(\nu_1+sN_1)(\nu_2+sN_2)} +
\frac{\chi(E_1\setminus C)}{\nu_1+sN_1}
+\frac{\chi(E_2\setminus C)}{\nu_2+sN_2} + \chi(X\setminus \{P\})  \\
&= \frac{-4}{(\frac15 - \frac15 s)(-s)} + \frac{-4}{\frac15 -
\frac15 s} + \frac 7{-s} + \chi(X\setminus \{P\}) = \frac {13}s +
\chi(X\setminus \{P\}) \, ,
\endalign
$$
yielding $\lim_{s\to 1} z_{st}(X;s)= z_{st}(X;1)=13 +
\chi(X\setminus \{P\})$.

\bigskip
\noindent {\bf 8.6.} {\sl Question.} Let $X$ be a $\Bbb
Q$-Gorenstein variety of arbitrary dimension without strictly log
canonical singularities. When is
$$\lim_{s \rightarrow 1} z_{st} (X;s) \in \Bbb Q \ ?$$

\bigskip
\bigskip
\noindent {\bf 9. Miscellaneous recent results}
\bigskip

Here we gather a collection of various results, which were
obtained after the redaction of the survey paper [DL8].
Undoubtedly some interesting work is missing, and this is of
course due to incompetence of the author of these notes. Any
suggestion is welcome.

\bigskip
\noindent $\bullet$ Aluffi [Al] noticed that the Euler
characteristic formula in (4.2) implies interesting similar
statements about Chern-Schwartz-MacPherson classes.

\smallskip
\noindent $\bullet$ Bittner [Bi2] calculated the relative dual of
the motivic nearby fibre and constructed a nearby cycle morphism
on the level of the Grothendieck group of varieties.

\smallskip
\noindent $\bullet$ More exotic motivic measures are introduced by
Bondal, Larsen and Lunts [BLL] and Drinfeld [Dr].

\smallskip \noindent $\bullet$  Using arc spaces and motivic
integration, Budur [Bu] relates the Hodge spectrum of a
hypersurface singularity to its {\sl jumping numbers} (which come
from multiplier ideals).

\smallskip \noindent $\bullet$ Campillo, Delgado and Gusein-Zade
[CDG1][CDG2][CDG3], and  Ebeling and Gusein-Zade [EG1][EG2]
studied filtrations on the ring of germs of functions on a germ of
a complex variety, defined by arcs on the singularity. An
important technique is integration with respect to the Euler
characteristic  over the projectivization of the space of function
germs; this notion is similar to (and inspired by) motivic
integration.

\smallskip \noindent $\bullet$
 Cluckers and Loeser [CL] built a more general theory for relative motivic
 integrals, avoiding moreover the completion of Grothendieck rings.
 These integrals specialize to both \lq classical\rq\ and arithmetic motivic integrals.

\smallskip \noindent $\bullet$ Dais and Roczen obtained formulas for
the stringy Euler number and stringy $E$-function for some special
classes of singularities [Da][DR].

\smallskip \noindent $\bullet$ Now available are the ICM 2002 survey
[DL9] and the recent  expository paper of Hales [Hal3] on the
theory of arithmetic motivic measure of Denef and Loeser [DL5].
Related work is in [DL10] and [Ni3].

\smallskip \noindent $\bullet$
 In [dSL] du Sautoy and Loeser associate motivic zeta functions
  to a large class of infinite dimensional Lie algebras.

\smallskip \noindent $\bullet$
 Ein,  Lazarsfeld, Musta\c t\v a and Yasuda have various other papers
 about spaces of jets, relating them for instance to singularities of
 pairs, in particular to the log canonical threshold, and to
 multiplier ideals [ELM][Mu2][Ya2].

\smallskip \noindent $\bullet$
 Koike and Parusi\'nski [KP] associated motivic zeta functions to real
 analytic function germs and showed that these  are  invariants of
 blow-analytic equivalence.  Fichou [Fi] obtained similar results in
 the context of Nash funcion germs. Both constructions are useful
 for classification issues.

\smallskip \noindent $\bullet$ Gordon [Go] introduced a motivic
analogue of the Haar measure for the (non locally compact) groups
$G(k((t)))$, where $G$ is a reductive algebraic groups, defined
over an algebraically closed field $k$ of characteristic zero.

\smallskip \noindent $\bullet$ Guibert [Gui] computed the motivic zeta
function  associated to irreducible plane curve germs, yielding a
new proof of the formula expressing the spectrum in terms of the
Puiseux data. Here he studied also a motivic zeta function for a
family of functions and related it with the Alexander invariants
of the family; this is used to obtain a formula for the Alexander
polynomial of a plane curve.

\smallskip \noindent $\bullet$ Guibert, Merle and Loeser [GML]
introduced iterated motivic vanishing cycles and proved  a motivic
 version of a conjecture of Steenbrink concerning the spectrum of
  hypersurface singularities.

\smallskip \noindent $\bullet$ Arithmetic motivic integration in the
context of $p$-adic orbital integrals and transfer factors is
considered by Gordon and Hales in [GH] and [Hal2]. An introduction
to this theory is [Hal1].

\smallskip \noindent $\bullet$ Ishii and Koll\'ar [IK] found counter
examples in dimensions at least 4 to the Nash problem, which
relates irreducible components of the space of arcs through a
singularity to exceptional components of a resolution. (And they
proved it in general for toric singularities.)

For a toric variety, Ishii [Is] described precisely the relation
between arc families and valuations, and obtained the answer to
the embedded version of the Nash problem.

\smallskip \noindent $\bullet$ Ito produced an alternative proof that
birational smooth minimal models have equal Hodge numbers [It1],
and that Batyrev's stringy $E$-function is well defined [It2],
using $p$-adic Hodge theory.

 \smallskip
\noindent $\bullet$ Kapranov [Ka] introduced another motivic zeta
function as the generating series for motivic measures of varying
$n$-fold symmetric products of a fixed variety. Larsen and Lunts
[LL1][LL2] determined for which surfaces this is a rational
function over $K_0(Var_{\Bbb C})$. It is not known whether it is
always a rational function over $\Cal M_{\Bbb C}$. See also [DL10,
\S7] and [BDN].

 \smallskip
\noindent $\bullet$
 For toric surfaces, Lejeune-Jalabert and Reguera [LR]
  and Nicaise [Ni1] computed an
explicit formula for the series $P(T)$ and $J(T)$, respectively.
This last paper also contains a sufficient condition for the
equality of $P(T)$ and the arithmetic Poincar\'e series of a toric
singularity, which is always satisfied in the surface case. A
counter
example for this equality in dimension 3 is given.

In [Ni2] Nicaise provides a concrete formula for $P(T)$ if the
variety has an embedded resolution of a simple form; this yields a
short proof of the formula for toric surfaces.

\smallskip \noindent $\bullet$ Loeser [Loe3] studied the behaviour of
motivic zeta functions of prehomogeneous vector spaces under
castling transformations; he deduced in particular how the motivic
Milnor fibre and the Hodge spectrum at the origin behave under
such transformations.

\smallskip \noindent $\bullet$
 Sebag [Se1][Se2] studied motivic integration and motivic zeta functions in the
 context of formal schemes. Loeser and Sebag [LS] developed a theory of motivic integration for smooth rigid
 varieties, obtained a motivic  Serre invariant, and  provided new
geometric birational invariants of degenerations of algebraic
varieties.

\smallskip \noindent $\bullet$
 Yasuda  [Ya1][Ya3] introduced {\sl twisted} jets and arcs over Deligne-Mumford
 stacks  and studied then motivic integration over them. As
 applications he obtained a McKay correspondence for general
 orbifolds (see also  [LP]), and a common generalization of the
 stringy $E$-function and the orbifold cohomology.

\smallskip \noindent $\bullet$
 Yokura [Yo] constructed Chern-Schwartz-MacPherson classes on
 pro-algebraic varieties and relates this to the motivic measure.

\enddocument